\documentclass[twoside,11pt,reqno]{amsart}

\usepackage[top=1.1in,bottom=1.1in,left=1.2in,right=1.2in]{geometry}
\usepackage{amsmath,amssymb,amsthm,amsfonts,enumerate, mathrsfs, mathtools, appendix}
\usepackage{relsize}
\usepackage{stmaryrd}
\usepackage{esint}
\usepackage{xcolor}
\usepackage{graphicx}
\usepackage{tikz} 
\usetikzlibrary{patterns,positioning,arrows,shapes,trees}

\usepackage[
colorlinks=true,
linkcolor=blue,
filecolor=blue,
urlcolor=red,
citecolor=magenta]{hyperref}

\usepackage{color}
\definecolor{MyDarkBlue}{RGB}{54,117,23}
\definecolor{MyDarkBlue}{cmyk}{0.8,0.3,0.8,0.4}
\definecolor{yellow}{rgb}{0.99,0.99,0.70}
\definecolor{white}{rgb}{1.0,1.0,1.0}
\definecolor{black}{rgb}{0.00,0.00,0.00}

\usepackage{fancyhdr}
{} 

\numberwithin{equation}{section}
\def\theequation{\arabic{section}.\arabic{equation}}
\newtheorem{theorem}{Theorem}[section]
\newtheorem{lemma}[theorem]{Lemma}
\newtheorem{remark}[theorem]{Remark}
\newtheorem{definition}[theorem]{Definition}
\newtheorem{proposition}[theorem]{Proposition}
\newtheorem{Examples}{Example}
\newtheorem{corollary}[theorem]{Corollary}
\newtheorem*{assumption}{Assumption}

\def\bt{\begin{theorem}}
\def\et{\end{theorem}}
\def\bl{\begin{lemma}}
\def\el{\end{lemma}}
\def\br{\begin{remark}}
\def\er{\end{remark}}
\def\bex{\begin{Examples}}
\def\eex{\end{Examples}}
\def\bd{\begin{definition}}
\def\ed{\end{definition}}
\def\bp{\begin{proposition}}
\def\ep{\end{proposition}}
\def\bc{\begin{corollary}}
\def\ec{\end{corollary}}
\def\ba{\begin{assumption}}
\def\ea{\end{assumption}}
\def\bpf{\begin{proof}}
\def\epf{\end{proof}}

\def\cH{{\mathcal H}}

\def\cL{{\mathcal L}}
\def\cM{{\mathcal M}}

\def\cP{{\mathcal P}}

\def\mD{{\mathbb D}}

\def\mK{{\mathbb K}}

\def\mN{{\mathbb N}}

\def\mR{{\mathbb R}}
\def\mS{{\mathbb S}}

\def\R{\mathbb R}
\def\N{\mathbb N}

\def\bP{{\mathbf P}}

\def\bP{{\mathbf P}}

\def\bE{{\mathbf E}}
\def\1{{\mathbf{1}}}

\def\sF{{\mathscr F}}

\def\sI{{\mathscr I}}

\def\eps{\varepsilon}
\def\d{\text{\rm{d}}}
\def\t{\tau}
\def\e{\mathrm{e}}

\def\a{\alpha}
\def\om{\omega}

\def\p{\partial}

\def\l{\lambda}
\def\si{\sigma}

\def\<{{\langle}}
\def\>{{\rangle}}

\def\dif{{\mathord{{\rm d}}}}

\def\={&\!\!=\!\!&}

\def\l{\left}
\def\r{\right}

\def\geq{\geqslant}
\def\leq{\leqslant}

\allowdisplaybreaks

\begin{document}
\title{Existence and Uniqueness for McKean-Vlasov equations with singular interactions}
\author{Guohuan Zhao}

\address{Institute of Applied Mathematics, Academy of Mathematics and Systems Science, CAS, Beijing, 100190, China}

\email{gzhao@amss.ac.cn}
	
\thanks{The research of Guohuan is supported by the National Natural Science Foundation of China (No. 12288201). }

\begin{abstract}
We investigate the well-posedness of following McKean-Vlasov equation in $\mathbb{R}^d$: 
$$
\mathrm{d} X_t=\sigma(t,X_t, \mu_{X_t})\mathrm{d} W_t+b(t, X_t, \mu_{X_t}) \mathrm{d} t,  
$$
where $\mu_{X_t}$ is the law of $X_t$. The existence of solutions is demonstrated when $\sigma$ satisfies certain non-degeneracy and continuity assumptions, and when $b$ meets some integrability conditions, and continuity requirements in the (generalized) total variation distance. Furthermore, uniqueness is established under additional continuity assumptions of a Lipschitz type. 


\noindent
  {\bf AMS 2020 Mathematics Subject: 60H10, 35K08, 35Q84} 
\end{abstract}

\maketitle
\section{Introduction}

\subsection{Motivations and Main results}
Let $\{W^i\}_{i\in \mN}$ be a collection of independent standard $d_1$-dimensional Wiener processes. Set $\Sigma: [0,1]\times\R^d \to \R^{d\times d_1}$, $B: [0,1]\times\R^d\to \R^d$. Consider the large systems of $N$ particles given by the coupled stochastic differential equations (SDEs)  
\begin{equation*}
\begin{aligned}
\d X^{i, N}_t =&\frac{1}{N-1}\sum_{j\neq i} \Sigma\left(t, X^{i, N}_t-X^{j, N}_t\right) \d W^{i}_t+\frac{1}{N-1}\sum_{j\neq i}B\left(t, X^{i, N}_t-X^{j, N}_t\right) \d t,  \\
X^{i, N}_0=&x^{i, N}_0 \ (i=1,2,\cdots N). 
\end{aligned}
\end{equation*}
The propagation of chaos phenomenon means that as $N$ goes to infinity, the limit of the empirical distribution $\frac{1}{N}\sum_{i=1}^N \delta_{X^{i,N}_t}$ is coincide with the marginal distribution of solution to the following McKean-Vlasov equation (MVE for short): 
\begin{equation}\label{Eq:MV}
    \d X_{t}=  \si(t,X_{t}, \mu_{X_t})\d W_t + b(t,X_{t}, \mu_{X_t}) \d t, \quad \mu_{X_0}=\pi, \tag{\bf{MV}} 
\end{equation}
where $W$ is a $d_1$-dimensional Brownian motion, $\mu_{X_t}$ is the law of $X_{t}$, and the maps 
\[
  \sigma: [0,1]\times\R^d\times \cP(\R^d) \to \R^{d\times d_1}~\mbox{ and }~ b:[0,1]\times\R^d\times \cP(\R^d) \to \R^{d}
\]
are given by 
\begin{equation}\label{Eq:A-B}
\si(t, x,m)=\int_{\R^d}\Sigma(t, x-y) m(\d y)\ \mbox{ and } \  b(t, x,m)=\int_{\R^d}B(t, x-y)m(\d y), 
\end{equation}
respectively (cf. \cite{sznitman1991topics}).

\medskip

In this paper, we conduct separate investigations into the existence and uniqueness issues concerning weak solutions to equations with a general form as described in \eqref{Eq:MV}, where the coefficients are not limited to the specific structure given by \eqref{Eq:A-B}. The formal definition of weak solutions is as follows 
\bd[weak solutions]\label{Def:weak}
Let $(\Omega, \sF, \{\sF_t\}_{t\in[0,1]}, \bP)$ be a filtered probability space satisfying common conditions and $(X,W)$ be a pair of adapted processes on it. We call $(\Omega, \sF, \{\sF_t\}_{t\in[0,1]}, \bP; X,W)$ is weak solution to \eqref{Eq:MV} if
\begin{enumerate}[(i)]
\item $\bP\circ X^{-1}_0=\pi$ and $W$ is a $d_1$-dimensional Brownian motion. 
\item For any $t\in [0,1]$, it holds that 
$$
X_t=X_0+\int_0^t \si(s, X_s, \mu_{X_s}) \d W_s+ \int_0^t b(s, X_s, \mu_{X_s}) \d s, \quad \bP-a.s.
$$
\end{enumerate}
\ed

In recent years, significant advancements have been made in the research on the well-posedness of \eqref{Eq:MV}. Here our attention is drawn to the literature \cite{rockner2021DDSDE}, where R\"ockner
and Zhang achieved remarkably comprehensive well-posedness results for \eqref{Eq:MV}, even in the presence of a singular drift term. However, it is worth noting that a limitation of their result is the requirement that the diffusion coefficient $\sigma$ should be independent of the measured variable $m$. One of the objectives of this article is to eliminate the limitation on $\sigma$ outlined in \cite{rockner2021DDSDE}. Additionally, we also attempt to extend previous results related to existence of solutions to \eqref{Eq:MV} to more general setting. Specifically, this work is focused on addressing the following two questions:
\begin{enumerate}[({\bf Q}$_1$)]
\item Is it possible to obtain  an existence result for \eqref{Eq:MV} when the drift term $b$ is singular in $t$ and $x$, and the coefficients $\sigma$ and $b$ are only continuous with respect to (w.r.t. in short) the measure variable $m$ in the {\bf total variation} distance? 
\item Can one extend the uniqueness result of R\"ockner-Zhang \cite[Theorem 4.3]{rockner2021DDSDE} for singular MVE to the case that $\sigma$ may also depend on the distribution variable $m$? 
\end{enumerate}

\medskip

Considering Peano's existence theorem for first order ordinary differential equations (ODEs in short), the first question is its natural analogy for distribution dependent SDEs. We employ the total variation distance rather than the widely used Wasserstein distance for a specific reason, as elucidated by the following observation: when the functions $\Sigma$ and $B$ in \eqref{Eq:A-B} are bounded, the coefficients $\sigma$ and $b$,  considered as functions of the third variable, may not exhibit continuity w.r.t. the measure variable $m$ in the Wasserstein metric. However, both $\sigma$ and $b$ are Lipschitz continuous w.r.t. $m$ in the total variation distance. Our answer to the first question is Theorem \ref{Thm:Ex} presented in Section \ref{Sec:Ex}, which indicates that much like Peano's existence results for ODEs, the presence of weak solutions for MVEs only requires the assumption that the coefficients is continuous in the (generalized) total variation distance.

As mentioned above, the response to the second question aims to supplement the uniqueness result presented in \cite{rockner2021DDSDE}. In fact, our motivation lies in the understanding that the scenario where the diffusion coefficient depends on the distribution is not just a theoretical extension; it has practical relevance. In fact, such a dependency of the diffusion coefficient on the distribution is prevalent in certain practical models, as evidenced in \cite{guilin2017kac} and its accompanying references. Our main result about uniqueness is Theorem \ref{Thm:Uni}, where we show that \cite[Theorem 4.3]{rockner2021DDSDE} can be extended to case that the map $m\mapsto \si(t,x,m)$  has a H\"older continuous linear functional derivatives (see section \ref{Sec:Pre} for the precise definition). 

\medskip

To give the reader a preliminary impression of the content of this article, we present one of our result below (see Theorem \ref{Thm:Ex} and Theorem \ref{Thm:Uni} for much more strong conclusions): 

\begin{proposition}\label{Prop:spf}
  Let  $\a\in (0,1)$, and $p,q\in (1,\infty)$ with $d/p+2/q<1$.  Assume $\si,b$ are given by \eqref{Eq:A-B}, $a=\frac{1}{2} \si\si^\top$ is uniformly elliptic, and 
$$
 \Sigma\in {L^\infty_tC^\a_x} ~\mbox{ and }~ B\in L^{q}_tL^p_x.  
$$
Then \eqref{Eq:MV} admits a unique weak solution. 
\end{proposition}

\subsection{Related literature}
The study of the general MVE (not limited to the special case of  \eqref{Eq:A-B}) has a long history, and there is a lot of literature focusing on the well-posedness problem. Among all, we mention that Funaki  \cite{funaki1984certain} proved the existence of martingale solutions to \eqref{Eq:MV} under some Lyapunov's type conditions as well as the uniqueness under global Lipschitz assumptions. By  Girsanov's transformation and Schauder's fixed point theorem,  Li-Min \cite{li2016weak} obtained the existence of weak solutions when $b$ is bounded and uniformly continuous in the Wasserstein distance w.r.t. to $m$, and the diffusion coefficient $\si$ depends only on $(t,x)$ is non-degenerate. Simultaneously, uniqueness was also proved in their work when $b$ is Lipschitz w.r.t.  the third variable (see also \cite{bauer2018strong}).  Under some one-side Lipschitz assumptions, Wang \cite{wang2018distribution} showed the strong well-posedness of \eqref{Eq:MV} and also some functional inequalities for the solutions. When $b$ only satisfies some integrability conditions and $\si,b$ are continuous in $m$, Huang-Wang \cite{huang2019distribution} proved the weak existence by using approximation argument and also strong uniqueness together with some standard conditions. It is crucial to emphasize that all of the mentioned results necessitate that $\sigma$ and $b$ exhibit at least uniform continuity concerning $m$ in the Wasserstein metric. 

As previously mentioned, in the context of  \eqref{Eq:A-B}, if $\Sigma$ and $B$ are merely bounded, the functions $\sigma$ and $b$, viewed in terms of the third variable, may not exhibit continuity w.r.t. $m$ in the Wasserstein distance. However, they are Lipschitz in the total variation distance. In this case, Shiga-Tanaka \cite{shiga1985central} proved the strong well-posedness for \eqref{Eq:MV} when $\si={\rm I}_{d\times d}$. Similar result was extended by Jourdain in \cite{jourdain1997diffusions}  for bounded drift $b$ with general form satisfying a Lipschitz assumption in $m$ w.r.t. the total variation metric. When the diffusion matrix $\si$ is uniformly non-degenerate and $\Sigma, B$ are  at most linear growth, Mishura-Veretennikov \cite{mishura2020existence} showed the existence of weak solutions. Meanwhile, they also proved the strong uniqueness by adding additional assumptions that $\si$ depends only on $t,x$ and is Lipschitz continuous in $x$. 
These results were  extended by Lacker in \cite{lacker2018strong} and later by R\"ockner-Zhang \cite{rockner2021DDSDE} to equations with possibly singular drifts in spaces $\cL^p_q$ (see Section \ref{Sec:Pre} for the definition). It is worth highlighting that these papers place significant emphasis on the Lipschitz-type requirement w.r.t. to measure variable $m$ for the drift coefficient, even when dealing with weak existence. Additionally, to facilitate the use of the Girsanov transformation for establishing uniqueness, these works necessitated assumptions of both the non-degeneracy of $\sigma$ and its independence from the distribution. 

The assumption that the diffusion coefficient is independent of the measure variable is relaxed by de Ranal and Frikha in \cite{de2020strong} and \cite{de2022well}. These studies contribute to the establishment of well-posedness results under specific conditions, utilizing a parametric matrix expansion of the transition density of the McKean-Vlasov process. The key conditions are as follows: 
\begin{enumerate}[(i)]
  \item The drift term $b$ is both bounded and Lipschitz continuous in $m$ when considering the total variation distance; 
  \item  The mapping $\sigma: [0,1]\times \mathbb{R}^d\times \mathcal{P}(\mathbb{R}^d)$ exhibits Hölder continuity in the linear functional derivatives concerning the measure variables. (For a precise definition, please refer to Section \ref{Sec:Pre}).
\end{enumerate}
The approach to proving uniqueness in this paper is significantly influenced by the insights derived from the two aforementioned articles. Additionally, in Section \ref{Sec:Uni}, we present a concrete example that underscores the necessity of continuity in the derivatives of linear functions.

\smallskip

We would also like to bring to the reader's attention some related papers on nonlinear Fokker-Planck equations, which bear a close connection to MVEs. Manita-Romanov-Shaposhnikov \cite{manita2015uniqueness} demonstrated the existence and uniqueness of solutions for the nonlinear Fokker-Planck equation \eqref{Eq:NFPE} under certain Lyapunov-type assumptions, employing a purely analytical argument. A similar endeavor was undertaken by Barbu and R\"ockner in \cite{barbu2018probabilistic, barbu2020nonlinear, barbu2021uniqueness}, where they considered nonlinear Fokker-Planck equations with coefficients depending on $m$ in a Nemytskii-type fashion. By applying Crandall-Liggett's theorem, they established some  existence results for these nonlinear PDEs. The subsequent weak existence of solutions to \eqref{Eq:MV} emerged as a result of applying the superposition principle, as elaborated in \cite{trevisan2012well} or \cite{figalli2008existence}.

\subsection{Main strategy}
Our approach is essentially rooted in the analysis of regularity and stability of heat kernels associated with the usual SDEs via Levi's parametrix method. Given two different second-order linear differential operators,  in Lemma \ref{Le:HKE},  we provide an estimate of the difference in their corresponding heat kernels, where the control term is given explicitly by the difference in the coefficients of the two operators. To attain the existence result, we consider the linearized SDE of \eqref{Eq:MV} given by 
$$
\d X_{t}^{\mu} = b(t, X_{t}^{\mu}, \mu_t) \d t + \si(t, X_{t}^{\mu}, \mu_t)  \d W_t, \quad \mathrm{law}(X_{0}^{\mu})=\pi , 
$$
where $\mu: (0,1]\to \cP(\R^d)$ is a curve in the probability space $\cP(\mR^d)$.  Subsequently,  by carefully choosing a suitable topological vector space $(V_\phi, d_{\phi})$ (see Section \ref{Sec:Ex}) and using the well-known Schauder's fixed point theorem, we show that the map $\psi: \mu\mapsto \{\mu_{X^\mu_t}\}_{t\in (0,1]}$ has at least one fixed point. This immediately yields our desired existence result (see Lemma \ref{Le:Ex} and Theorem \ref{Thm:Ex} below).  It must be pointed out that the above approach to prove the weak existence is  mainly inspired by \cite{ chen2017heat}, \cite{li2016weak} and \cite{zhang2020stochastic}. 

For uniqueness, before giving our answer to question ({\bf Q}$_2$), we first present an example to show that the uniqueness result might fail,  even if the diffusion coefficient $a=\frac{1}{2}\si\si^\top$ is uniformly elliptic and its linear functional derivative (see section \ref{Sec:Pre}) of $\sigma$ is bounded (which implies the Lipschitz continuity of $\sigma$ in $m$ w.r.t. the total variation distance).  To establish uniqueness, we draw inspiration from \cite{de2022well}, introducing additional H\"older regularity assumptions on the linear functional derivative of the diffusion coefficient $\si$. Subsequently, again leveraging Levi's parametrix expansion, the desired result is proved by estimating the difference of two solutions' transition probability densities. 

\subsection{ Organization of the article}

In Section \ref{Sec:Pre}, we review some basic facts that will be used in this article.  Section \ref{Sec:Ex} and Sections \ref{Sec:Uni} contain the proof of Theorem \ref{Thm:Ex} (Existence) and Theorem \ref{Thm:Uni} (Uniqueness) respectively. For the sake of completeness and consistency, we provide an appendix that contains proofs of the weak well-posedness of SDEs with irregular coefficients (Lemma \ref{Le:LSDE}) and certain properties of the Kato functions (Lemma \ref{Le:Kato}).

We closed this section by collecting  some frequently used notations. 
\begin{itemize}
\item The letter $C$ denotes a constant, whose value may change in different places.
\item We use  $A\lesssim B$ and $A\asymp B$ 
to denote $A\leq C B$ and $C^{-1} B\leq A\leq CB$ for some unimportant constant $C>0$, respectively.
\item Suppose $x$ is a vector in a Euclidean space and $A$ is a matrix, 
$$|x|:=\Big(\sum_{i} |x_i|^2\Big)^{1/2}, \quad  |A|:=\max_{i,j} |A_{ij}|.
$$
\item $B_R:=\l\{x\in\R^d: |x|<R\r\}$,  $\mD:=\l\{(s,x, t,y): 0\leq s<t\leq 1, x,y\in \R^d\r\}$. 
\item Give $\Lambda>1$, $S_\Lambda$ is the collection of $d\times d$ symmetric matrices whose eigenvalue are between $\Lambda^{-1}$ and $\Lambda$. For any $\a(0,1)$ and $N>0$, define 
$$
  \mS(\Lambda, \a, N):= \l\{a: [0,1]\times \R^d\to S_\Lambda: \|a_{ij}\|_{L^\infty([0,1]; C^\a(\R^d))}\leq N, \ i,j=1,\cdots, d \r\}. 
$$
\item for any $f:A\to \mR$ and $g: I\times A \to \mR$, set 
\[
  \|f\|_{L^p(A)} := \l(\int_{A} |f(x)|^p \d x\r)^{1/p}, \quad \|g\|_{L^{q,p}(I\times A)}:= \l[ \int_{I}\l(\int_{A} |g(t,x)|^{p} \d x\r)^{q/p}\r]^{1/q}. 
\]
\item $\chi\in C_c^\infty(B_2)$, $\chi\in [0,1]$ and $\chi\equiv1$ in $B_1$;  $\chi_z(x):=\chi(x-z)$. 
\item For $\beta\in [0,2)$, we introduce the index set $\sI_\beta$ as following:
$$
\sI_\beta := \l\{(p,q): p,q \in [2,\infty),  \ \tfrac{d}{p}+\tfrac{2}{q}<2-\beta\r\}. 
$$
\item $\phi: \R^d\to [1,\infty)$ is a smooth, radial and increasing function. 
\item 
$\cM(\R^d)$ ($\cP(\R^d)$) is the collection of signed (probability) measures on $\R^d$. 
\item The Wasserstein distance $W_p$ for $p\geq 1$ is defined as 
$$
W_p(m, m') := \inf_{\pi\in \Pi(m,m')} \l(\int_{\R^d\times \R^d} |x-y|^p \pi(\d x, dy)\r)^{1/p}, 
$$
where $\Pi (m ,m')$ denotes the collection of all measures on $\R^d\times \R^d$ with marginals $m$  and $m'$  on the first and second factors respectively. 
\end{itemize}

\section{Some preparations}\label{Sec:Pre}
In this section, we make some
preparations. We first introduce the localized Bessel potential spaces for later use. 

For $s\geq 0$ and $p\in (1,\infty)$, the usual Bessel potential space $H^{s, p}$ is defined as 
$$
H^{s, p}=\l\{g\in L^p: ({\rm I}-\Delta)^{s/2} g\in L^p\r\}. 
$$
Recall that $\chi\in C_c^\infty(B_2)$, $\chi\in [0,1]$ and $\chi\equiv1$ in $B_1$, and $\chi_z(x):=\chi(x-z)$. Suppose that $g\in L^p_{loc}(\R^d)$ and $f\in L^{q,p}_{loc}(\R^{d+1})$, we denote 
\begin{align*}
\|g\|_{\cL^{p}}:= \sup_{z\in \R^d} \|g\chi_z\|_{L^p(\R^d)}
~\mbox{ and }~ \|f\|_{\cL^{p}_{q} (T)}:= \sup_{z\in \R^d} \|f \chi_z\|_{L^q([0,T]; L^p(\R^d))}. 
\end{align*}
We also introduce the localized $H^{s, p}$-space: 
$$
\cH^{s, p}:= \Big\{g\in H^{s, p}_{loc}: \| g \|_{\cH^{s, p}}:=\sup_{z\in \R^d} \|g\chi_z\|_{H^{s,p}}<\infty \Big\}
$$
and the localized space-time function space $\cH^{s,p}_q(T)$ with norm 
$$
\|f\|_{\cH^{s, p}_q(T)}:= \sup_{z\in \R^d} \l(\int_0^T \|f(t)\chi_z\|_{H^{s,p}}^q \d t\r)^{1/q}. 
$$
When $T=1$, $\|f\|_{\cL^{p}_q(1)}$ and $\|f\|_{\cH^{s, p}_{q}(1)}$ are denoted by $\|f\|_{\cL^{p}_{q}}$ and $\|f\|_{\cH^{s, p}_{q}}$, respectively.  
 
\medskip
 
Two main ingredients of proving the existence of solutions to \eqref{Eq:MV} (see Lemma \ref{Le:Ex} and Proposition \ref{Prop-Key}) are the Fr\'echet-Kolmogorov theorem and the Schauder-Tychonoff fixed point theorem, which are presented below.

\bl[Fr\'echet-Kolmogorov  theorem]\label{Le:FK}
Let $K$ be a bounded set in $L^p(\R^d)$ with $p\in [1,\infty)$. The subset $K$ is relatively compact if and only if the following properties hold: 
\begin{enumerate}[(i)]
\item $\lim_{R\to \infty} \sup_{f\in K} \|f\1_{B^c_R}\|_{L^p} =0. $
\item 
$\lim_{h\to 0} \sup_{f\in K}\|f(\cdot+h)-f\|_{L^p}=0$. 
\end{enumerate}
\el

\bl[Schauder-Tychonoff fixed point theorem, \cite{pata2019fixed}]\label{Le:ST}
Let $V$ be a locally convex Hausdorff topological vector space, $S$ be a nonempty closed convex subset of $V$, $\psi$ be a continuous mapping on $S$. If  $K=\psi(S)$ is a relatively compact subset of $S$, then $\psi$ has a fixed point in $K$.
\el
Next, we give a lemma about the weak well-posedness of usual SDEs, which will be used to prove the weak existence of solutions to MVE. 

\bl\label{Le:LSDE}
Let 
$$
  \sigma: [0,1]\times \R^d\to \R^{d\times d_1}, \quad b: [0,1]\times \R^d\to \R^d. 
$$
Let $(p,q)\in \sI_1$, $\alpha\in (0,1)$ and $ \Lambda, N_1, N_2>1$. Assume that 
$a=\frac{1}{2}\si\si^\top\in \mS(\Lambda, \a, N_1)$ and $\| b\|_{\cL^p_q}\leq N_2$, then equation 
\begin{equation}\label{Eq:SDE1}
X_{s,t} = X_{s,s}+\int_s^t b(r, X_{s,r})\d r + \int_s^t \si(r, X_{s,r}) \d W_r, \quad 0\leq s\leq t\leq 1
\end{equation}
has a unique weak solution. 
\el

In order to address the heat kernels of diffusion processes with singular drifts, we employ certain generalized Kato's function spaces, which were initially introduced in \cite{chen2017heat} and \cite{zhang2018singular}. Let $I$ represent an interval within $\R_+$, and consider a measurable function $f: I\times \R^d\to \R$. Noting that throughout this discussion and in the subsequent context, we consistently extend the function $f$ to $\R^{d+1}$ by setting $f(t,x)=0$ when $t$ is outside the interval $I$.  For any $\beta\geq 0, \lambda>0$, define 
$$
\eta_\beta (t,x):= (\sqrt{t}+|x|)^{-d-\beta}, \quad t> 0, x\in \R^d
$$
and 
\begin{align*}
K^\beta_{f}(T):=&\sup_{(t,x)\in\mR^{d+1}}\int_0^T\!\!\!\int_{\R^d} \eta_\beta(s, y) |f(t+s, x+y)|\d y\dif s\\
&+\sup_{(t,x)\in\mR^{d+1}}\int_0^T\!\!\!\int_{\R^d} \eta_\beta(s, y)|f(t-s, x-y)| \d y\dif s, \quad T>0. 
\end{align*}
The generalized Kato's class is defined by 
$$
\mK^\beta:=\left\{f: \mR^{d+1}\to \R \ \big|  \quad \lim_{\delta\to 0}K^\beta_{f}(\delta)=0\right\}.
$$
For any $\lambda>0, \gamma\in \R$, put  
$$
\varrho_{\lambda, \gamma}(t, x):= t^{(-d+\gamma)/2} \e^{-\lambda |x|^2/t} , \quad t> 0, x\in \R^d. 
$$
$\varrho_{\lambda, 0}$ is denoted by $\varrho_\lambda$ for simplicity.  

\medskip 

The following facts will be used frequently. 
\bl\label{Le:Kato}
\begin{enumerate}[(i)]
\item For any $\gamma\geq 0$ and $\kappa\in (0,1)$, 
\begin{align*}
    |x|^\gamma \varrho_{\lambda, 0}(t, x) \lesssim \varrho_{\kappa\lambda, \gamma}(t, x). 
\end{align*}
\item For any $\lambda>0$, $\beta\geq 0$ and $0<t\leq 1$, $x\in \R^d$, 
\begin{equation}\label{Eq:rhoVeta}
\varrho_{\lambda, -\beta}(t,x)\lesssim \eta_\beta(t,x). 
\end{equation}
\item Let $\beta\in (0,2)$ and $(p,q)\in \sI_\beta$. Then for any $f\in \cL^p_q$ and $T\in [0,1]$, 
\begin{equation}\label{Eq:KvsL}
K_{f}^\beta(T)\lesssim T^{\frac{1}{2}(2-\beta-\frac{d}{p}-\frac{2}{q})}  \|f\|_{\cL^p_q(T)}. 
\end{equation} 
\item Let $\beta\geq \beta'\geq 0$. For any $0\leq s<t\leq 1$ and $x,y\in \R^d$, 
\begin{equation}\label{Eq:kato}
\begin{aligned}
&\int_{s}^{t}\!\!\!\int_{\mathbb{R}^{d}} \varrho_{\lambda,-\beta'}(\tau-s, x-z)  |b(\tau, z)| \varrho_{2 \lambda,-\beta}(t-\tau, z-y) \d z\d \tau \\
\lesssim&  K_{|b|}^\beta(t-s) \varrho_{\lambda,-\beta'}(t-s, x-y). 
\end{aligned}
\end{equation}
\end{enumerate}
\el

\medskip

Finally, we introduce the concept of derivatives w.r.t. the measure variable, a notion commonly employed in linearization procedures. It is worth recalling that $\phi: \R^d\to [1,\infty)$ is a smooth, radial, and increasing function. Giving  $f : \cP_\phi(\R^d) \to \R$, we say $f$ has a linear functional derivative if there exists a function $\frac{\delta f}{\delta m}: \cP_\phi(\R^d) \times \R^d \to \R$ such that for any compact subset $K$ of $\cP_\phi(\R^d)$, 
 $$
 \sup_{m\in K}\frac{\delta f}{\delta m}(m)(y)\leq C_K\phi(y), \quad \forall y\in \R^d 
 $$ 
 and for any $m, m'\in \cP_\phi(\R^d)$,
\begin{align*}
\lim _{\varepsilon \downarrow 0} \frac{f\left((1-\varepsilon) m+\varepsilon m^{\prime}\right)-f(m)}{\varepsilon}=\int_{\mathbb{R}^{d}}\frac{\delta f}{\delta m}(m)(y) \left(m^{\prime}-m\right)(\d y). 
\end{align*}
Note that $\frac{\delta f}{\delta m}$ is defined up to an additive constant. We adopt the normalization convention 
$$
\int_{\R^d} \frac{\delta f}{\delta m}(m)(y) m(\d y)=0.
$$ 
For any $m, m' \in \mathcal{P}_{\phi}(\mathbb{R}^{d})$, we have 
\begin{equation}\label{Eq:f1-f2}
f(m)-f(m')=\int_{0}^{1} \!\!\!\int_{\mathbb{R}^{d}} \frac{\delta f}{\delta m}\left(\lambda m+(1-\lambda) m'\right)(y) \left(m-m'\right)(\d y) \d \lambda, 
\end{equation}
which implies 
\begin{equation*}
    |f(m)-f(m')|\leq \sup_{\lambda\in[0,1]; y\in \R^d} \l |\frac{\delta f}{\delta m}\left(\lambda m+(1-\lambda) m'\right)(y) \r| \|m-m'\|_{{\rm TV}}. 
\end{equation*}
See also \cite[section 2.2]{cardaliaguet2015master} for more details. 

\section{Weak existence}\label{Sec:Ex}
In this section, following the introduction of some definitions, we present our primary result concerning the existence of weak solutions to equation \eqref{Eq:MV} in Lemma \ref{Le:Ex}. To attain this objective, we also offer an intriguing technical result in Lemma \ref{Le:HKE}. This lemma addresses the regularity and stability characteristics of heat kernels associated with second-order operators featuring singular first-order terms.
  
\medskip
  
Recall that $\phi: \R^d\to [1,\infty)$  is a smooth, radial and increasing function. Here and below, we further assume that for each $\lambda>0$, 
\begin{equation}\label{Eq:growth}
\sup_{t\in[0,1]; h\in B_1}\int_{\R^d} (|\phi|+|\nabla \phi|)(x-h-y) \varrho_\lambda(t,y)\d y\lesssim_\lambda \phi(x),  \tag{{\bf G}}
\end{equation}

Two typical examples of functions satisfying \eqref{Eq:growth} are 
\[
  \phi(x)=\exp({\sqrt{1+|x|^2}}) ~ \mbox{ and }~ \phi(x)=(1+|x|^2)^p\, (p\geq 0).
\]

For any $m\in \cM(\R^d)$, define 
$$
\<f, m\>:= \int_{\R^d} f(x)m(\d x), \quad \|m\|_{\phi}:=\<\phi, |m|\>, 
$$
where $|m|$ is the variation of $m$. Let 
$$
\cM_\phi(\R^d):=\l\{m\in \cM(\R^d): \|m\|_{\phi}<\infty\r\} ~ \mbox{ and }~\cP_\phi(\R^d):= \cM_\phi(\R^d)\cap\cP(\R^d). 
$$

\br\label{Rek-metric}
\begin{enumerate}
\item Obviously, $\|m-m'\|_{{\rm TV}}\leq \|m-m'\|_\phi$; 
\item  If $\phi(x)=(1+|x|^2)^{1/2}$, then by \cite[Theorem 6.15]{villani2008optimal}, it holds that $W_1(m, m')\leq \|m-m'\|_{\phi}. $
\end{enumerate}
\er

Now we can state our result on existence. 
\bt[Existence]\label{Thm:Ex}
Let $\a\in (0,1)$, $p,q\in (1,\infty)$ with $d/p+2/q<1$. Assume $\om(\delta): [0,2]\to \R_+$ is an increasing function with $\lim_{\delta\downarrow0}\om(\delta)= 0$ and  $\ell(t,\delta): [0,1]\times [0,2] \to \R_+$ is another nonnegative function, which is increasing in $\delta$, and satisfies 
\begin{equation}\label{eq:mod}
\|\ell(\cdot, 2)\|_{L^q([0,1])}<\infty ~\mbox{ and }~ \lim_{\delta\downarrow 0}\|\ell(\cdot, \delta)\|_{L^q([0,1])}= 0. 
\end{equation}
Suppose that 
$a=\frac{1}{2}\si\si^\top$ is uniformly elliptic, 
\begin{align*}
    \sup_{(t,m)\in [0,1]\times \cP(\R^d)}  |\sigma(t,x,m)-\sigma(t,x',m)|\leq C |x-x'|^\a,
\end{align*}
\begin{equation}\label{Eq:Esig}
    \sup_{(t,x)\in [0,1]\times\R^d} |\sigma(t,x,m)-\sigma(t,x,m')|\leq \om (\|m-m'\|_{\phi}) \tag{{\bf E}$_\sigma$}
\end{equation}
and that 
\begin{equation}\label{Eq:Eb}
    \|b(\cdot,\cdot, \delta_0)\|_{\cL^p_q}<\infty, \quad \|b(t,\cdot, m)-b(t, \cdot, m')\|_{\cL^p} \leq \ell(t, \|m-m'\|_{\phi}). \tag{{\bf E}$_b$}
\end{equation}
Then for any $\pi\in \cP_{\phi}(\mR^d)$, equation \eqref{Eq:MV} has at least one weak solution. Here $\|\|_{\phi}$  (see Section \ref{Sec:Ex} for the definition of $\|\cdot\|_{\phi}$ and $\cP_{\phi}(\mR^d)$). 
\et

Before presenting the proof of Theorem \ref{Thm:Ex}, we begin by introducing a vector space $V_\phi$. This vector space encompasses all continuous mappings from the half-open interval $(0,1]$ to $\cM_\phi(\R^d)$ and includes one of its closed convex subsets, denoted as $S_\phi$. The latter subset will assume a pivotal role in our proof.
\bd 
Set 
\begin{align*} 
V_\phi:=&C((0,1]; \cM_\phi(\R^d))=\l\{\mu:(0,1]\to \cM_\phi(\R^d)\Big| \ \lim_{t\to t_{0}} \|\mu_t-\mu_{t_0}\|_\phi=0, \ \forall t_0\in (0,1]\r\}.
\end{align*}
The distance $d_\phi$ on $V_\phi$ is defined as 
\begin{align*}
d_\phi (\mu, \mu'):= \max_{k\in \N_+} \frac{2^{-k}\sup_{t\in[\frac{1}{k},1]}\|\mu_t-\mu'_t\|_\phi }{(1+\sup_{t\in [\frac{1}{k},1]}\|\mu_t-\mu'_t\|_\phi )}. 
\end{align*}
Put 
\begin{align*}
S_\phi:=&\l\{\mu\in V_\phi: \mu_t\in \cP_{\phi}(\R^d) \mbox{ and } \mu_t(\d x)\ll \d x\mbox{ for each }t\in(0,1]\r\}.  
\end{align*}  
\ed
By definition, one can see that 
\begin{enumerate}
\item A sequence $\{\mu^n\}\subseteq V_\phi$ converges to $\mu$ in $V_\phi$ iff 
$$
\lim_{n\to\infty} \sup_{t\in[t_0, 1]} \|\mu_t^n-\mu_t\|_\phi= 0,\quad \forall t_0\in (0,1]; 
$$
\item $S_\phi$ is a convex and closed subset of $V_\phi$. 
\end{enumerate}

\begin{remark}
  The definitions of $V_\phi$ and $S_\phi$ may appear non-standard since we exclude the endpoint '0'. There are two key reasons for this choice: 
  \begin{enumerate}
    \item  Even the nice curve $ \mu: [0,1]\ni t\mapsto \mu_{W_t}\in \cP(\R^d)$ is not continuous at  $0$ if $\cP(\R^d)$ is equipped with total variation norm; 
    \item As mentioned in the introduction, our existence result will be established through the application of the Schauder-Tychonoff fixed-point theorem. Therefore, we require a suitable topological space in which all compact sets can be distinctly identified. Characterizing the compact sets within $V_\phi$ is a straightforward task.
  \end{enumerate}
\end{remark}

\medskip

Given $\mu\in S_\phi$, put 
$$
\si^\mu(t,x):= \si(t,x,\mu_t),\quad  a^\mu=\frac{1}{2}\si^{\mu} (\si^\mu)^t, \quad  b^\mu(t,x):= b(t,x,\mu_t). 
$$
The main technical result of this section is 
\begin{lemma}\label{Le:Ex}
Let $(p,q)\in \sI_1$, $\gamma_0=1-\frac{d}{p}-\frac{2}{q}>0$, $r\in (2/\gamma_0, \infty]$, $\alpha\in (0,1)$, $\Lambda>1$ and $N_1, N_2>1$. 
Suppose for each $\mu\in S_\phi$, 
\begin{align}\label{Aspt1}
a^\mu\in \mS(\Lambda, \a, N_1), \quad \| b^\mu\|_{\cL^p_q}\leq N_2 {\tag {\bf E$'_1$}}
\end{align}
and 
\begin{align}\label{Aspt2}
\mu\mapsto (a^\mu, b^\mu)\mbox{ is  continuous from }(S_{\phi}, d_\phi)  \mbox{ to } \cL^\infty_r\times \cL^p_q. {\tag {\bf E$'_2$}} 
\end{align}
Then for any $\pi\in \cP_\phi(\R^d)$, equation \eqref{Eq:MV} has at least one weak solution.
\end{lemma}

Using Lemma \ref{Le:Ex}, we can present  
\bpf[Proof of Theorem \ref{Thm:Ex}]
We only need to verify that $a,b$ satisfy \eqref{Aspt1} and \eqref{Aspt2}.  
By the assumptions on $\sigma$, one sees that $a^\mu\in \mS(\Lambda, \a, N_1)$ for some $N_1>0$, and \eqref{Eq:Eb} implies that for each $\mu\in S_\phi$ and $t_0\in (0,1]$ 
\begin{equation}\label{eq:b-lplq}
\begin{aligned}
\|b^{\mu}\|_{L^q([0,t_0]; \cL^p)} \leq& \|b(t,x,\mu_t)-b(t, x, \delta_0)\|_{L^q_t([0,t_0]; \cL^p_x)} + \|b(\cdot, \delta_0)\|_{L^q([0,t_0]; \cL^p)}\\
\leq& \|\ell(t,2)\|_{L^q_t([0,t_0])}+ \|b(\cdot, \delta_0)\|_{L^q([0,t_0]; \cL^p)}. 
\end{aligned}
\end{equation}
Thus, $a, b$ satisfies \eqref{Aspt1}. 

Assume that $\mu^n \to \mu$ in $S_\phi$, then for each $t_0\in (0,1]$, $
\sup_{t\in [t_0,1]}\| \mu^n_t -\mu_t\|_{\phi} \to 0 \, (n\to \infty)$.  \eqref{Eq:Esig} implies that for any $r\in (1,\infty)$
\begin{align*}
\|a^{\mu^n}-a^{\mu}\|_{\cL^\infty_r}\leq   2 N_1 t_0^{1/r}  + \| a(t,x, \mu^n_t)-a(t,x,\mu_t)\|_{L^\infty_{t,x}([t_0,1]\times \R^d)}. 
\end{align*}
Letting $n\to \infty$ and then $t_0\to 0$ in the right side of above inequality, we obtain $\|a^{\mu^n}-a^\mu\|_{\cL^\infty_r} \to 0 \ (n\to \infty)$. Similarly, 
\begin{align*}
\|b^{\mu_n}-b^\mu\|_{\cL^p_q}\leq& \|b^{\mu_n}\|_{L^q([0, t_0]; L^p)} +\|b^\mu\|_{L^q([0, t_0]; L^p)}+\|b^{\mu_n}-b^\mu\|_{L^q([t_0, 1]; L^p)} \\
\overset{\eqref{eq:b-lplq}}{\leq} & 2  \l(\|\ell(t,2)\|_{L^q_t([0,t_0])}+ \|b(\cdot, \delta_0)\|_{L^q([0,t_0];  \cL^p)} \r)+\Big\| \|b(t,\cdot,\mu^n_t)-b(t,\cdot,\mu_t)\|_{\cL^p} \Big\|_{L^q_t([t_0,1])} \\
\leq&2 \|\ell(t,2)\|_{L^q_t([0,t_0])}+ 2 \|b(\cdot, \delta_0)\|_{L^q([0,t_0];  \cL^p)} + \| \ell(t, \|\mu^n_t-\mu_t\|_{\phi}) \|_{L^q_t([t_0,1])}. 
\end{align*}
Letting $n\to\infty$ and then $t_0\to 0$, by \eqref{eq:mod} and \eqref{Eq:Eb}, we obtain $\|b^{\mu_n}-b^\mu\|_{\cL^p_q}\to0 \ (n\to\infty)$. Therefore, $a, b$ also fulfil \eqref{Aspt2}. So, we complete our proof. 
\epf

The strategy of proof for Lemma \ref{Le:Ex} is following: 
\begin{enumerate}
    \item 
    According to Lemma \ref{Le:LSDE}, for any $\mu\in S_{\phi}$, equation 
    \begin{equation}\label{Eq:SDE2}
      \d X_{t}^{\mu}= b^\mu(t, X_{t}^{\mu}) \d t + \si^\mu(t, X_{t}^{\mu})  \d W_t, \quad \mu_{X_0}=\pi\in \cP_{\phi}(\mR^d) 
    \end{equation} 
    admits a unique weak solution $X^\mu$, provided that $\si$ and $b$ satisfy \eqref{Aspt1}. Set 
    $$
    \psi: S_\phi\ni \mu\mapsto \{\mu_{{X^\mu_{t}}}\}_{t\in (0,1]}; 
    $$ 
    \item 
    Verify that the triple $(V_\phi, S_\phi, \psi)$ satisfies all the conditions in Lemma \ref{Le:ST}, so that one can find a fixed point of $\psi$ in $S_\phi$.
\end{enumerate}
Therefore, the main task is to show 
\bp\label{Prop-Key}
Under the same conditions of Lemma \ref{Le:Ex}. Assume that $\psi$ is the map defined above and $K:=\psi(S_\phi)$, then 
\begin{enumerate}
\item  $K\subseteq S_\phi\subseteq V_\phi$;  
\item   $K$ is relatively compact in $S_\phi$;  
\item the map $\psi: S_{\phi}\to K$ is continuous. 
\end{enumerate}
\ep

Utilizing Propostion \ref{Prop-Key}, it is easy to show 
\bpf[Proof of Lemma \ref{Le:Ex}]
By  \cite[Theorem 1.37 and Remark 1.38]{rudin1973functional}, $(V_\phi,d_\phi )$ is a locally convex topological vector space, and obviously, $S_\phi$ is a closed convex subset of $V_\phi$. Thanks to Schauder-Tychonoff fixed point theorem and Proposition \ref{Prop-Key}, $\psi$ has a fixed point $ \mu=\psi( \mu)\in K$. By our assumptions and Lemma \ref{Le:LSDE}, we deduce the desired result. 
\epf

In order to establish our Proposition \ref{Prop-Key}, we require some auxiliary lemmas. The following lemma offers a basic estimation of the difference between the determinants of two positive-definite matrices. 
\bl
Given $\Lambda \geq 1$, $a\geq b>0$. For any $d\times d$ symmetric positive-definite matrices $A_a, \widetilde A_a, B_b$ satisfying 
$$
a \Lambda^{-1} {\rm I}\leq A_a, \widetilde A_a\leq a \Lambda {\rm I}, \quad b \Lambda^{-1} {\rm I}\leq B_b\leq b \Lambda {\rm I}, 
$$
it holds that 
\begin{equation}\label{Eq:detA-detB}
|\det A_a-\det \widetilde A_a| \leq C  a^{d-1} |A_a-\widetilde A_a|
\end{equation}
and 
\begin{equation}\label{Eq:detA+B-detA}
|\det (A_a+B_b)-\det A_a |\leq C a^{d-1} b, 
\end{equation}
where $C$ only depends on $d$ and $\Lambda$. 
\el
\bpf
We only prove \eqref{Eq:detA-detB} here, since the proof for \eqref{Eq:detA+B-detA} is similar. We can assume $\delta:= |A_a-\widetilde A_a|\leq \frac{a}{100d\Lambda}$, otherwise \eqref{Eq:detA-detB} is obviously true. Suppose $QA_aQ^t=D$, where $Q$ is an orthogonal matrix and $D={\rm diag}(a\lambda_1,\cdots, a\lambda_d)$ with $\lambda_i\in [\Lambda^{-1}, \Lambda]$. Then 
$$
Q\widetilde A_a Q^t= Q(\widetilde A_a- A_a) Q^t + {\rm diag} (a\lambda_1, \cdots, a \lambda_d), 
$$
and 
$$
\det A_a = \Pi_{k=1}^d a\lambda_k, \quad 
\det \widetilde A_a = \Pi_{k=1}^d a\lambda_{k}+ p(\eps_{ij}), 
$$
where $\eps_{ij}=[ Q(\widetilde A_a- A_a) Q^t ]_{ij}$ and $p$ is a polynomial of $\eps_{ij}$ without zero order term. Note that  
\begin{align*}
|\eps_{ij}|=&\sum_{k, l} |Q_{ik}|\, |(\widetilde A_a- A_a)_{kl}| \, | Q^t_{l, j}|\\
\leq&  \frac{\delta}{2}\sum_{k, l} (|Q_{ik}|^2+|Q_{lj}|^2 ) \leq \delta(d+1)/2\leq d\delta, 
\end{align*}
 we obtain 
\begin{align*}
\l| \det \widetilde A_a -\det A_a\r| =
|p(\eps_{ij})|\leq C(d)(a\Lambda+d\delta)^{d-1} (d\delta) \leq C(\Lambda, d)a^{d-1}\delta. 
\end{align*}
\epf

The next lemma is about the regularity and stability properties of heat kernels associated with second order elliptic operators with singular first order terms. 
\bl\label{Le:HKE}
Let $(p,q)\in \sI_1$,  $\gamma_0:=1-\frac{d}{p}-\frac{2}{q}>0$, $\alpha\in (0,1)$ and $ \Lambda, N_1, N_2>1$. Assume   $a=\frac{1}{2}\si\si^\top\in \mS(\Lambda, \a, N_1)$ and $\| b\|_{\cL^p_q}\leq N_2$. Then for each $x\in \R^d$, \eqref{Eq:SDE1}  admits a unique (weak) solution $X_{s,t}(x)$ $(X_{s,s}=x)$ and the law of $X_{s,t}(x)$ has a density $p(s, x; t, y)$. Moreover,
\begin{enumerate}
\item there is a constant $\lambda\in (0,1)$ depending only on $d,\a, p, q, \Lambda, N_1, N_2$ such that 
\begin{enumerate}[(i)]
\item {\em (Gaussian estimate):}  for all $0\leq s<t\leq 1$ and $x,y\in\mR^d$,
\begin{align}\label{Eq:TSE}
\begin{aligned}
\varrho_{\lambda^{-1}}(t-s, x-y)\lesssim p(s, x; t, y)\lesssim  \varrho_{\lambda}(t-s, x-y); 
\end{aligned}
\end{align}
\item {\em (H\"older estimate in $t$ and $y$):} for any $\gamma \in (0, \a\wedge \gamma_0)$, $0\leq s<t_1<t_2\leq 1$ and $x,y, y_1, y_2\in\mR^d$, it holds that 
\begin{align}\label{Eq:Holder-t}
\begin{aligned}
|p(s,x; t_2, y)-p(s, x; t_1, y)|\leq C |t_1-t_2|^{\frac{\gamma}{2}} \sum_{i=1}^2 \varrho_{\lambda, -\gamma}(t_i-s, x-y)\end{aligned}
\end{align}
and 
\begin{align}\label{Eq:Holder-y}
\begin{aligned}
|p(s,x; t, y_1)-p(s, x; t, y_2)|
\leq C |y_1-y_2|^{\gamma} \sum_{i=1}^2\varrho_{\lambda, -\gamma}(t-s, x-y_i), 
\end{aligned}
\end{align}
where $C$ only depends on $d, \a,  p, q, \Lambda, N_1, N_2$ and  $\gamma$. 
\end{enumerate}
\item Assume $\widetilde a\in \mS(\Lambda, \a, N_1)$, $\|\widetilde b\|_{\cL^p_q}\leq N_2$ and $\widetilde p$ is the heat kernel associated with  $\widetilde L=\widetilde a_{ij}\p_{ij}+\widetilde b_i\p_i$. Then for any $r\in (2/\gamma_0, \infty]$ and $\eta\in (\frac{2}{2+\a r},1)$, there is a constant $C$ only depends on $d,\a, p, q, r, \eta, \Lambda, N_1, N_2$ such that for all $0\leq s<t\leq 1$ and $x, y\in \R^d$, 
\begin{align}\label{Eq:stable}
|p-\widetilde p|(s,x;t,y) \lesssim \l(\| a-\widetilde a\|_{\cL^\infty_r}^{1-\eta} +\|b-\widetilde b\|_{\cL^p_q}\r)  \varrho_{\lambda, -2/r}(t-s, x-y).
\end{align}
\end{enumerate}
\el
\bpf
We only present the proof of \eqref{Eq:Holder-t} and\eqref{Eq:stable} below, since \eqref{Eq:TSE} is already proved in \cite{chen2017heat} and the proof for \eqref{Eq:Holder-y} is similar with \eqref{Eq:Holder-t}. 

We use the classic Levi's parametrix method to prove our conclusions, and point out that it is enough to prove the result when $0\leq t-s \leq T$ is small. This is because we can use the reproducing property of the fundamental solution to cover the case of  $0\leq s <t\leq 1$. 

(1) 
For any $p^{(1)}, p^{(2)}, \cdots, p^{(n)}: \mD\to \R$, set
\begin{align*}
&\l[p^{(1)}\otimes p^{(2)}\otimes \cdots \otimes p^{(n)}\r](s,x; t,y)\\
:=& \int_{s<\tau_1<\cdots<\tau_{n-1}<t}\int_{\R^{nd}}  p^{(1)}(s,x; \tau_1, z_1)p^{(2)}(\tau_1,z_1; \tau_2, z_2)\\
&\qquad \qquad \qquad \quad \cdots p^{(n)}(\tau_{n-1}, z_{n-1}; t, y) \d z_1\cdots\d z_{n_1}~\d \tau_1\cdots\d \tau_{n-1}. 
\end{align*}

For any fixed $y\in \R^d$, put 
\[
  L_0f(s, x) := a_{ij}(s, y)\p_{x_i x_j}f(s,x). 
\]
Set $A_{s,t}(y):=\int_s^t a(\tau, y) \d \tau$ and 
$$
p_0(s,x; t, y):=\frac{\mathrm{e}^{-\left\langle A_{s,t}^{-1}(y) (x-y), (x-y)\right\rangle }}{\sqrt{(4 \pi)^{d} \operatorname{det}(A_{s,t}(y))}}. 
$$
Then $p_0$ satisfies  
$$
(\p_s+L_0)p_0=0, \quad  x\in \R^d,\ a.e.\, s\in (0,t]
$$
(cf. \cite{friedman2008partial}). Since the generator of $X$, denoted by $L=a_{ij}\p_{ij}+b_i\p_i$ can be viewed as a perturbation of $L_0$ by $L-L_0$, heuristically the fundamental solution (or heat kernel) $p(s, x; t, y)$ of $L$ should satisfy the following Duhamel’s formula: 
\[
  p= p_0+p\otimes (L-L_0)p_0. 
\]
Formally, we have 
\begin{align}\label{Eq:p}
p=\sum_{n=0}^\infty p_n :=\sum_{n=0}^\infty p_0 \otimes [(L-L_0)p_0]^{\otimes^n},  
\end{align}
Rigorously,  our assumption about $b$, in conjunction with Lemma \ref{Le:Kato} (iii), implies that $b$ belongs to $\mK^{1+\gamma}$ for some $\gamma \in [0, \gamma_0)$. Based on this and the discussion in the Appendix of \cite{chen2017heat}, it is confirmed that Equation \eqref{Eq:p} indeed holds. To demonstrate \eqref{Eq:Holder-t}, as outlined in \cite{chen2017heat}, we rely on \eqref{Eq:p} and carefully estimate $|p_n(s, x; t_2, y) - p_n(s, x; t_1, y)|$. 

{\bf Claim $1$}: for all $n\in \N$, $0\leq s<t\leq s+T$ and $x,y\in \R^d$, 
\begin{equation}\label{eq:pn}
|p_n(s, x ; t, y)| \leq \lambda_n \varrho_{\lambda}(t-s; x-y), 
\end{equation}
where 
$$
\lambda_n:= C_1^n \l(T^{\frac{\a}{2}}+K_{|b|}^1(T)\r)^n, 
$$
and $C_1=C_1(d, \a, p, q, \Lambda, N_1, N_2)>0$ and $T\in (0,1]$ are constants that will be determined later. 
\eqref{eq:pn} can be established through induction as follows: basic calculation yields that there is a constant $\lambda >0$ only depending on $d, \Lambda$ such that for any $k\in \{0,1,2\}$
\begin{align}\label{Eq:kp0}
|\nabla_x^k p_0(s,x;t,y)|  \lesssim  \varrho_{3\lambda, -k}(t-s, x-y). 
\end{align}
By the H\"older regularity of $a(t, \cdot)$ and \eqref{Eq:kp0}, we have 
\begin{equation}\label{Eq:L-L0}
|(L-L_0)p_0|(s,x;t,y)\lesssim \varrho_{\lambda, \a-2}(t-s, x-y)+ |b(s,x)| \varrho_{2\lambda, -1}(t-s, x-y). 
\end{equation}
Assume that \eqref{eq:pn} is already proved for some $n\geq 0$, by the above estimate, one then can see that 
\begin{align*}
|p_{n+1}|(s,x;t,y)=& |p_n\otimes (L-L_0)p_0|(s,x;t,y)\\
\overset{\eqref{eq:pn},\eqref{Eq:L-L0}}{\leq}& \lambda_n  \int_{s}^{t}\!\!\!\int_{\R^d}\varrho_{\lambda}(\tau-s, x-z)  \varrho_{\lambda, \a-2}(\tau,z;t,y) \d z\d \tau\\
&+\lambda_n  \int_{s}^{t}\!\!\!\int_{\R^d}\varrho_{\lambda}(\tau-s, x-z)\l|b(\tau,z)\r| \varrho_{2\lambda, -1}(t-\tau, z-y) \d z\d \tau\\
\overset{\eqref{Eq:kato}}{\lesssim} & \lambda_n \l( \int_s^t (t-\tau)^{\frac{\a}{2}-1} \d \tau+ K_{|b|}^1(t-s)\r) \varrho_{\lambda}(t-s, x-y) \\
\lesssim &  \lambda_n  (T^{\frac{\a}{2}}+ K_{|b|}^1(T)) \varrho_{\lambda}(t-s, x-y),
\end{align*}
i.e. there is a constant $C$ only depends on $d, \a, \Lambda, N_1$ such that 
$$
|p_{n+1}|\leq \lambda_n C \l(T^{\frac{\a}{2}}+K_{|b|}^1(T)\r) \varrho_{\lambda}= \lambda_{n+1} \varrho_{\lambda},    
$$
thus \eqref{eq:pn} holds for all $n\in \mN$. 

\smallskip

{\bf Claim $2$}: for any $\gamma\in (0, \a\wedge \gamma_0)$, $n\in \N$, $0\leq s<t_1< t_2\leq T$ and $x,y\in \R^d$, it holds that 
\begin{equation}\label{Eq:pn-pn}
\begin{aligned}
&|p_{n+1}(s,x;t_1,y)-p_{n+1}(s,x;t_2,y)| \\
\lesssim&\lambda_n|t_1-t_2|^{\frac{\gamma}{2}} \sum_{i=1}^2 \varrho_{2\lambda, -\gamma}(t_i-s, x-y). 
\end{aligned}
\end{equation}
We first prove \eqref{Eq:pn-pn} for $n=0$. It is easy to see that one only need to consider the case of $0\leq s<t_1<t_2\leq 1\wedge (2t_1-s)$.  By the elementary inequality $|\e^{x}-\e^{y}|\leq |x-y| (\e^{x}+\e^{y})$, \eqref{Eq:detA+B-detA} and noticing that $t_2-s\leq 2(t_1-s)$, we have 
\begin{equation*}
\begin{aligned}
&|p_0(s,x; t_1, y)-p_0(s,x; t_2, y)| \\
\lesssim & \l|(\det A_{s,t_1}(y))^{-1/2}-(\det A_{s,t_2}(y))^{-1/2}\r| \e^{-\<A^{-1}_{s,t_1}(y)(x-y), x-y\>}\\
&+ \l|(\det A_{s,t_2})^{-1/2}\r| \l| \e^{-\<A^{-1}_{s,t_1}(y)(x-y), x-y\>}-\e^{-\<A^{-1}_{s,t_2}(y)(x-y), x-y\>}\r|\\
\lesssim &\e^{-\<A^{-1}_{s,t_1}(y)(x-y), x-y\>}   |\det A_{s,t_1}(y)-\det A_{s,t_2}(y)|\cdot\\
&\l[\sqrt{\det A_{s, t_1}(y)\det A_{s, t_2}(y)} \l(\sqrt{\det A_{s, t_1}(y)}+\sqrt{\det A_{s, t_2}(y)}\r)\r]^{-1} \\
&+ (t_2-s)^{-d/2} \l|A_{s,t_1}^{-1}(y)[A_{s,t_2}(y)-A_{s,t_1}(y)]A_{s,t_2}^{-1}(y)\r| |x-y|^2 \\
&\quad \quad\quad\quad \quad \quad \quad \quad\quad \quad \quad \quad \l| \e^{-\<A^{-1}_{s,t_1}(y)(x-y), x-y\>}+\e^{-\<A^{-1}_{s,t_2}(y)(x-y), x-y\>}\r|\\
\overset{\eqref{Eq:detA+B-detA}}{\lesssim}&  |t_2-t_1| (t_1-s)^{-1}\varrho_{3\lambda, 0} (t_1-s, x-y)+  |t_2-t_1| \frac{|x-y|^2}{(t_2-s)^2} \varrho_{3\lambda, 0} (t_2-s, x-y)\\
\lesssim & |t_1-t_2|^{\frac{\gamma}{2}} \sum_{i=1}^2 \varrho_{2\lambda, -\gamma}(t_i-s, x-y).
\end{aligned}
\end{equation*}
Similarly, for all $0\leq s<t_1\leq t_2 \leq 1$ and $k=0,1,2$, one can show that 
\begin{equation}\label{eq:p0t1-p0t2}
\begin{split}
&|\nabla^k p_0(s,x; t_1, y)-\nabla^k p_0(s,x; t_2, y)| \\
\lesssim & |t_1-t_2|^{\frac{\gamma}{2}}\l[\varrho_{2\lambda, -k-\gamma}(t_1-s, x-y) + \varrho_{2\lambda, -k-\gamma}(t_2-s, x-y) \r]. 
\end{split}
\end{equation}
Now assume that \eqref{eq:pn-pn} is true for some $n\geq 0$, then by definition 
\begin{align*}
\begin{aligned}
&|p_{n+1}(s,x;t_1,y)-p_{n+1}(s,x;t_2,y)| \\
=&| p_n\otimes (L-L_0)p_0(s,x,t_1,y)-p_n\otimes (L-L_0)p_0(s,x,t_2,y)| \\
\lesssim & \int_s^{t_1}\!\!\!\int_{\R^d} |p_n(s,x;\tau,z)|\cdot | a_{ij}(\tau, z)-a_{ij}(\tau, y)|\\
&\quad \quad \quad \quad \quad \quad \quad \quad \quad \quad  \cdot |\p_{x_ix_j} p_0(\tau, z;t_1, y)-\p_{x_ix_j} p_0(\tau, z;t_2, y)| \d z\d\tau \\
&+ \int_{t_1}^{t_2}\!\!\!\int_{\R^d} |p_n(s,x;\tau,z)|\cdot |a_{ij}(\tau, z)-a_{ij}(\tau, y)|\cdot |\p_{x_ix_j} p_0(\tau, z;t_2, y)|\d z\d\tau\\
&+\int_s^{t_1}\!\!\!\int_{\R^d} |p_n(s,x;\tau,z)|\cdot |b_i(\tau, z)| \cdot |\p_{x_i} p_0(\tau,z;t_1,y)-\p_{x_i} p_0(\tau,z;t_2,y)| \d z\d \tau \\
&+\int_{t_1}^{t_2}\!\!\!\int_{\R^d} |p_n(s,x;\tau,z)|\cdot |b_i(\tau,z)| \cdot |\p_{x_i} p_0(\tau,z; t_2, y)| \d z\d\tau\\ 
=:& I_1^n+I_2^n+I_3^n+I_4^n. 
\end{aligned}
\end{align*}
Using \eqref{eq:pn} and \eqref{eq:p0t1-p0t2}, we have
\begin{equation*}
\begin{aligned}
 I_1^n\lesssim& \lambda_n  |t_1-t_2|^{\frac{\gamma}{2}} \int_s^{t_1}\!\!\!\int_{\R^d} \varrho_\lambda(\tau-s, x-z)\sum_{i=1}^2\varrho_{\lambda, \a-\gamma-2}(t_i-\tau, z-y)\d z\d \tau \\
\lesssim& \lambda_n  |t_1-t_2|^{\frac{\gamma}{2}}\l[\int_s^{t_1} (|t_1-\tau|^{\frac{\a-\gamma}{2}-1}+|t_2-\tau|^{\frac{\a-\gamma}{2}-1}) \d \tau\r]   \sum_{i=1}^2\varrho_{\lambda}(t_i-s, x-y)\\
\lesssim&  \lambda_n  |t_1-t_2|^{\frac{\gamma}{2}} T^{\frac{\a-\gamma}{2}} \sum_{i=1}^2\varrho_{\lambda}(t_i-s, x-y). 
\end{aligned}
\end{equation*}
Similarly, 
\begin{align*}
I_2^n\overset{\eqref{eq:pn},\eqref{Eq:kp0}}{\lesssim}&\lambda_n \int_{t_1}^{t_2} |t_2-\tau|^{\frac{\a}{2}-1}\d \tau \sum_{i=1}^2\varrho_{\lambda}(t_i-s, x-y) \\
\lesssim& \lambda_n   |t_1-t_2|^{\frac{\gamma}{2}} T^{\frac{\a-\gamma}{2}}\sum_{i=1}^2\varrho_{\lambda} (t_i-s; x-y)
\end{align*}
and 
\begin{align*}
I_3^n\overset{\eqref{eq:pn},\eqref{eq:p0t1-p0t2}}{\lesssim}
&\lambda_n  |t_1-t_2|^{\frac{\gamma}{2}}  \int_{s}^{t_1}\!\!\!\int_{\R^d} \varrho_{\lambda}(\tau-s,x-z) |b(\tau,z)|  \sum_{i=1}^2\varrho_{2\lambda, -1-\gamma} (t_i-\tau; z-y)\d z\d \tau \\
\overset{\eqref{Eq:kato}}{\lesssim} &\lambda_n  |t_1-t_2|^{\frac{\gamma}{2}}  K_{|b|}^{1+\gamma}(T)\sum_{i=1}^2\varrho_{\lambda} (t_i-s; x-y),  
\end{align*} 
and 
\begin{align*}
I_4^n\overset{\eqref{Eq:kp0},\eqref{eq:pn}}{\lesssim} &\lambda_n|t_2-t_1|^{\frac{\gamma}{2}}   \int_{t_1}^{t_2}\!\!\! \int_{\R^d} \varrho_{\lambda}(\tau-s,x-z)  |b(\tau,z)|  \varrho_{2\lambda, -1-\gamma}(t_2-\tau, z-y) \d z\d\tau \\
\overset{\eqref{Eq:kato}}{\lesssim}  &  \lambda_n |t_1-t_2|^{\frac{\gamma}{2}}   K_{|b|}^{1+\gamma}(T) \sum_{i=1}^2\varrho_{\lambda} (t_i-s; x-y). 
\end{align*}
Combining all the estimates above, we get \eqref{Eq:pn-pn}. 

Having \eqref{Eq:pn-pn} at our disposal, we can observe that
\begin{align*}
&|p(s,x;t_1, y)-p(s,x;t_2, y)|\leq \sum_{n=0}^\infty |p_n(s,x;t_1, y)-p_n(s,x;t_2, y)|\\
\overset{\eqref{Eq:pn-pn}}{\lesssim}&  |t_1-t_2|^{\frac{\gamma}{2}}  \l[\sum_{i=1}^2\varrho_{\lambda,-\gamma} (t_i-s; x-y)+
\sum_{n=0}^\infty \lambda_n\sum_{i=1}^2\varrho_{\lambda} (t_i-s; x-y) \r]. 
\end{align*}
Recall that $\lambda_n= C_1^n \l(T^{\frac{\a}{2}}+K_{|b|}^1(T)\r)^n$, by choosing $T$ sufficiently small so that $\sum_{n=0}^\infty \lambda_n <\infty$, we obtain \eqref{Eq:Holder-t}. 

(2). Now we start to prove \eqref{Eq:stable}. Again we use the expansion \eqref{Eq:p} and estimate each $|p_n-\widetilde p_n|$ separately. 

Denote 
\[
  \delta(t):= \|a(t,x)-\widetilde a (t,x)\|_{L^\infty_x} ~ \mbox{ and }~ \delta_a:=\| a-\widetilde a\|_{L^r([0,1];L^\infty)}.
\]
By definition, 
\begin{align}\label{eq:p0-tildep0}
\begin{aligned}
&|p_0(s,x;t,y)-\widetilde p_0(s,x;t,y)| \\\lesssim& \l| (\det A_{s,t}(y))^{-1/2}-(\det \widetilde A_{s,t}(y))^{-1/2}\r| \e^{- \<A_{s,t}^{-1}(y)(x-y), x-y\>} \\
&+ (\det \widetilde A_{s,t}(y))^{-1/2} \l|\e^{- \<A_{s,t}^{-1}(y)(x-y), x-y\>}-\e^{- \<\widetilde A_{s,t}^{-1}(y)(x-y), x-y\>}\r|. 
\end{aligned}
\end{align}
 Like the proof for \eqref{eq:p0t1-p0t2}, using  \eqref{Eq:detA-detB}, \eqref{eq:p0-tildep0} and H\"older's inequality, we have 
\begin{align}
&|p_0(s,x;t,y)-\widetilde p_0(s,x;t,y)| \notag \\
\lesssim &\l[(t-s)^{-d} |\det A_{s,t}(y)-\det \widetilde A_{s,t}(y)|+ |A^{-1}_{s,t}(y)-\widetilde A^{-1}_{s,t}(y) | |x-y|^2\r] \varrho_{3\lambda, 0}(t-s, x-y) \notag \\
\overset{\eqref{Eq:detA-detB}}{\lesssim}& \l[(t-s)^{-1} +\frac{|x-y|^2}{(t-s)^2} \r] \l|A_{s,t}(y)-\widetilde A_{s,t}(y)\r| \varrho_{3\lambda, 0}(t-s,x-y) \notag \\
\lesssim& \int_s^t |a(\tau, y)-\widetilde a(\tau, y)|\d \tau \cdot \varrho_{3\lambda, -1}(t-s, x-y) \label{eq:p0-p0-1} \\
\lesssim& \delta_a (t-s)^{1-\frac{1}{r}}  \varrho_{3\lambda, -1}(t-s, x-y)\lesssim  \delta_a \varrho_{2\lambda, -2/r}(t-s, x-y). \notag
\end{align}
Similarly, for any $k=0,1,2$, we have 
\begin{align}\label{Eq:2p0-p0}
\begin{aligned}
&|\nabla_x^k p_0(s,x;t,y)-\nabla_x^k \widetilde p_0(s,x;t,y)| \lesssim \delta_a \varrho_{2\lambda, -k-2/r}(t-s, x-y). \end{aligned}
\end{align}

{\bf Claim $3$}: for each $n\in \mN$, $0\leq s<t\leq T$ and $x,y\in \R^d$, 
\begin{equation}\label{eq:pn-pn}
|p_n-\widetilde p_n|(s, x ; t, y) \leq   \kappa_n(T) (\delta_a^{1-\eta} +\delta_b)  \varrho_{\lambda , -2/r}(t-s; x-y), 
\end{equation}
where 
$$
\delta_a=\| a-\widetilde a\|_{L^r_tL^\infty_x}, \quad \delta_b:=\|b-\widetilde b\|_{\cL^p_q},\quad  \eta\in (\tfrac{2}{2+\a r},1)
$$
and $\kappa_n(T)\, (\forall n\in \mN)$ is a constant that will be determined later. 

Again, we prove \eqref{eq:pn-pn} by induction. Assume \eqref{eq:pn-pn} holds for some $n\in \mN$, by \eqref{Eq:p}, we have 
\begin{align*}
&p_{n+1}-\widetilde p_{n+1}= p_n\otimes (L-L_0)p_0-\widetilde  p_n\otimes (\widetilde L-\widetilde L_0)\widetilde p_0\\
=& (p_n-\widetilde p_n) \otimes (L-L_0)p_0+ \widetilde p_n \otimes [(L- L_0)-(\widetilde L-\widetilde L_0)]p_0+ \widetilde  p_n\otimes (\widetilde L-\widetilde L_0)(p_0- \widetilde p_0)\\
=&: J_1+J_2+J_3. 
\end{align*}
Noting that $\frac{d}{p}+\frac{2}{q}=1-\gamma_0<2-(1+\frac{2}{r})$, by \eqref{Eq:KvsL}, we can find a continuous function $\om: [0,1]\to \R_+$  depending only on $d, p, q, r$ such that  
\begin{equation}\label{eq:kvsl}
K^{1+\frac{2}{r}}_{|b|}(t)+K^{1+\frac{2}{r}}_{|\widetilde b|}(t)\leq N_2 \om(t), \quad K_{|b-\widetilde b|}^1(t)\leq \|b-\widetilde b\|_{\cL^p_q}\ \om (t) 
\end{equation}
and $\om(t)\to 0$ as $t\to 0$.  So for $J_1$, by the fact that $\frac{\a}{2}-\frac{1}{r}>0$, one sees that for all $0\leq s\leq t\leq T$ and $x,y\in \mR^d$, 
\begin{align*}
|J_1|(s,x; t,y)  \overset{\eqref{Eq:L-L0}, \eqref{eq:pn-pn}}{\lesssim}& \kappa_n(T)(\delta_a^{1-\eta} +\delta_b)   \int_{s}^{t}\!\!\!\int_{\R^d} \varrho_{\lambda, -2/r}(\tau-s, x-z) \varrho_{\lambda, \a-2}(t-\tau, z-y) \d z \d \tau \\
&+  \kappa_n(T) (\delta_a^{1-\eta} +\delta_b)     \int_{s}^{t}\!\!\!\int_{\R^d} \varrho_{\lambda, -2/r}(\tau-s, x-z) |b(\tau,z)| \varrho_{2\lambda, -1}(t-\tau, z-y) \d z \d \tau \\
\overset{\eqref{Eq:kato}}{\lesssim}  &   \kappa_n(T) (\delta_a^{1-\eta} +\delta_b)   \l( \int_s^t (\tau-s)^{-\frac{1}{r}}(t-\tau)^{\frac{\a}{2}-1} \d \tau \r) \varrho_{\lambda}(t-s, x-y) \\
&+   \kappa_n(T) (\delta_a^{1-\eta} +\delta_b)   K_{|b|}^1(T) \varrho_{\lambda, -2/r}(t-s, x-y)\\
\overset{\eqref{eq:kvsl}}{\lesssim}&   \kappa_n(T) (\delta_a^{1-\eta} +\delta_b)   \l( T^{\frac{\a}{2}-\frac{1}{r}} \varrho_{\lambda}(t-s, x-y)+  \om(T) \varrho_{\lambda, -2/r}(t-s, x-y) \r), 
\end{align*}

For $J_2$, recalling that $\delta(t)= \|a(t,x)-\widetilde a (t,x)\|_{L^\infty_x}$, by our assumptions on $a$ and $\widetilde a$,  it holds that 
\begin{align*}
&\l| (a_{ij}-\widetilde a_{ij})(t,x)-(a_{ij}-\widetilde a_{ij})(t, y)\r| \\
\leq&  2  \delta(t)  \wedge 2N_1 |x-y|^\a  \lesssim \delta^{1-\eta}(t) |x-y|^{\a \eta},  
\end{align*}
thus
\begin{align*}
&\l | [(L- L_0)-(\widetilde L-\widetilde L_0)]p_0 \r| (s,x;t,y)\\
\lesssim& \delta^{1-\eta}(s) \varrho_{\lambda, \a\eta -2} (t-s, x-y) + |b-\widetilde b|(s,x) \varrho_{2\lambda, -1}(t-s, x-y). 
\end{align*}
Using the above estimate and H\"older's inequality, we obtain 
\begin{align*}
|J_2|(s,x; t,y) \lesssim& \kappa_n(T) \int_s^t \delta^{1-\eta}(\tau) \d \tau \int_{\R^d} \varrho_{\lambda}(\tau-s, x-z)  \varrho_{\lambda, \a\eta-2}(t-\tau,z-y) \d z \\
&+ \kappa_n(T)\int_{s}^{t}\!\!\!\int_{\R^d} \varrho_{\lambda}(\tau-s, y-z) |b-\widetilde b|(\tau,z) \varrho_{2\lambda, -1}(t-\tau,z-y)\d z \d \tau\\
\overset{\eqref{Eq:kato}}{\lesssim}& \kappa_n(T) \l[ \int_s^t \delta^{1-\eta}(\tau)(t-\tau)^{\frac{\a\eta}{2}-1} \d \tau + K^1_{|b-\widetilde b|}(T)\r] \varrho_{\lambda}(t-s,x-y)\\
\lesssim& \kappa_n(T) \l(\int_s^t |\delta(\tau)|^r \d \tau \r)^{\frac{1-\eta}{r}} \l(\int_s^t (t-\tau)^{\frac{r(\a\eta-2)}{2(r+\eta-1)}} \d \tau \r)^{\frac{r+\eta-1}{r}} \varrho_{\lambda}(t-s,x-y)\\
&+ \kappa_n(T) K_{|b-\widetilde b|}^1(T)\varrho_{\lambda}(t-s,x-y)\\
\overset{\eqref{eq:kvsl}}{\lesssim} & \kappa_n(T)  \l(\delta_a^{1-\eta}  T^{\frac{\a\eta}{2}-\frac{1-\eta}{r}}+  \delta_b\ \om(T) \r) \varrho_{\lambda}(t-s,x-y), 
\end{align*}
where we use the fact $\frac{r(\a\eta-2)}{2(r+\eta-1)}>-1$. 

For $J_3$, by the $\a$-H\"older regularity of $x\mapsto a(t,x)$ and \eqref{Eq:2p0-p0}, we get 
\begin{align*}
&(\widetilde L-\widetilde L_0)(p_0-\widetilde p_0) (s,x;t,y)\\
\lesssim & \delta_a \varrho_{\lambda, \a-2/r-2}(t-s, x-y)+ \delta_a  |\tilde b(s,x)|\varrho_{2\lambda, -1-2/r}(t-s, x-y). 
\end{align*}
Hence, 
\begin{align*}
|J_3|(s,x;t,y)\lesssim & \kappa_n(T)\delta_a   \int_s^t (t-\tau)^{\frac{\a}{2}-\frac{1}{r}-1}\d \tau \int_{\R^d} \varrho_{\lambda}(\tau-s, x-z)\varrho_{\lambda}(t-\tau, z-y) \d z\\
&+ \kappa_n(T)\delta_a   \int_s^t \!\!\!\int_{\R^d}\varrho_{\lambda}(\tau-s, x-z) |\widetilde b(\tau, z) |\varrho_{2\lambda, -1-2/r}(t-\tau, z-y) \d z \d \tau\\
\overset{\eqref{Eq:kato},\eqref{eq:kvsl}}{\lesssim}  & \kappa_n(T)\delta_a    \l(T^{\frac{\a}{2}-\frac{1}{r}}+ \om (T) \r)\varrho_{\lambda}(t-s, x-y). \end{align*}
Combining the above estimates, we can see that there is a constant $C_2$ depends only on $d, \a, p, q, r, \eta, \Lambda, N_1, N_2$ such that 
\begin{align*}
\begin{aligned}
|p_{n+1}-\widetilde p_{n+1}| \leq &C_2 \l(T^{\frac{\a\eta}{2}-\frac{1-\eta}{r}}+ \om(T) \r) \kappa_n(T)(\delta_a^{1-\eta} +\delta_b)    \varrho_{\lambda, -2/r}. 
\end{aligned}
\end{align*}
If we define
\begin{equation}\label{eq:kappa}
    \kappa_n(T):= C_2^n\l(T^{\frac{\a\eta}{2}-\frac{1-\eta}{r}}+ \om(T) \r)^n, \quad \forall n\in\mN
\end{equation}
where $C_2=C_2(d, \alpha, p, q, r, \eta, \Lambda, N_1, N_2)>0$, and $T\in (0,1]$ to be determined later, then for all $0\leq s\leq t\leq T$ and $x,y\in \mR^d$, we have
\begin{align*}
|p_{n+1}-\widetilde p_{n+1}|(s,x;t,y) \leq \kappa_{n+1}(T)(\delta_a^{1-\eta} +\delta_b) \varrho_{\lambda, -2/r}(t-s,x-y).
\end{align*}
So \eqref{eq:pn-pn} holds for all $n\in \mN$ if $\kappa_n(T)$ given by \eqref{eq:kappa}. Choosing $T$ sufficiently small such that $C_2 \l(T^{\frac{\a\eta}{2}-\frac{1-\eta}{r}}+ \om(T) \r) <1$, we obtain 
\begin{align*}
|p-\widetilde p| \leq \sum_{n=0}^\infty |p_n-\widetilde p_n|\lesssim\sum_{n=0}^\infty \kappa_n(T) (\delta_a^{1-\eta} +\delta_b)  \varrho_{\lambda, -2/r}\lesssim (\delta_a^{1-\eta} +\delta_b)  \varrho_{\lambda, -2/r}. 
\end{align*}
\epf

Now, we are in a position to provide the proof for the Proposition \ref{Prop-Key}. 
\bpf[Proof of Proposition \ref{Prop-Key}]
(1). Given $\mu\in S_\phi$, let $p^\mu(s, x;t, y)$ be the heat kernel associated with $L^\mu:=a^\mu_{ij}\p_{ij}+b^\mu_i\p_i$. Recall that  $\psi(\mu)_t={\rm law}(X_{t}^\mu)$ is the one dimensional distribution of the unique weak solution to \eqref{Eq:SDE2}.  For any $t\in [0,1]$, 
\begin{align*}
\<\phi, \psi(\mu)_t\>=&\bE \phi(X^\mu_t)
\lesssim \int_{\R^d} \pi(\d x) \int_{\R^d} \phi(y) p^\mu(0,x;t,y) \d y\\
\overset{\eqref{Eq:TSE}}{\lesssim} &\int_{\R^d} \pi(\d x) \int_{\R^d} \phi(y) \varrho_\lambda (t,x-y) \d y \overset{\eqref{Eq:growth}}{\lesssim} \<\phi, \pi\> <\infty, 
\end{align*}
which implies $\psi(\mu)_t\in \cP_\phi(\R^d)$. Let $\gamma_0=1-\frac{d}{p}-\frac{2}{q}>0$ and $\gamma \in (0, \a\wedge \gamma_0)$. For any $0<t_0\leq t_1\leq t_2\leq 1$ and $f\in L^\infty(\R^d)$,  by the Markov property,  
\begin{align*} 
&\l|\int_{\R^d} f(y)\phi(y) [\psi(\mu)_{t_2}-\psi(\mu)_{t_1}](\d y)\r|=| \bE f\phi(X^\mu_{t_2})-\bE f\phi(X^{\mu}_{t_1}) | \\
=&\l| \int_{\R^d} \pi(\d x)\int_{\R^d}[p^\mu(0,x; t_2, y)- p^\mu(0,x;t_1, y)](f\phi)(y) \d y\r| \\
\leq& \|f\|_\infty \int_{\R^d} \pi(\d x)\int_{\R^d}\l| p^\mu(0,x; t_2, y)- p^\mu(0,x;t_1, y)\r|\phi(y) \d y \\
\overset{\eqref{Eq:Holder-t}}{\lesssim}&\|f\|_\infty  |t_1-t_2|^{\frac{\gamma}{2}} \sum_{i=1}^2 \int_{\R^d} \pi(\d x) \int_{\R^d}   \phi(y)  \varrho_{\lambda, -\gamma}(t_i,x-y) \d y\\
\overset{\eqref{Eq:growth}}{\lesssim}_{\gamma, t_0} &  \|f\|_\infty t_0^{-\frac{\gamma}{2}} |t_1-t_2|^{\frac{\gamma}{2}}, 
\end{align*}
which implies 
\begin{equation}\label{eq:equicont}
\|\psi(\mu)_{t_2}-\psi(\mu)_{t_1}\|_{\phi} \lesssim_{\gamma, t_0}  |t_1-t_2|^{\frac{\gamma}{2}}. 
\end{equation}
Thus, $K=\psi(S_{\phi})\subseteq S_\phi$. For each $t_0\in (0,1]$ and $\mu\in S_\phi $, denote  
$$
\mu|_{t_0}: [t_0,1]\ni t\mapsto \mu_t\in \cP(\R^d), \quad K|_{t_0}:= \l\{\mu|_{t_0}: \mu\in K\r\}. 
$$
\eqref{eq:equicont} also implies the equicontinuity of $K|_{t_0}$.  

\medskip

(2). By the definition of $(V_\phi, d_\phi )$ and the standard  diagonal argument, we only need to show that $K|_{t_0}$ is relatively compact in $C([t_0, 1]; \cM_\phi(\R^d))$, which equipped with norm $\|\mu|_{t_0}\|= \sup_{t\in [t_0,1]} \|\mu_t\|_\phi $. By Arzela-Ascoli's Theorem and \eqref{eq:equicont}, we only need to prove that for each fixed $t\in (0, 1]$, $\{\psi(\mu)_t: \mu\in S_\phi\}$ is a relatively compact set in $(\cP_\phi(\R^d), \|\cdot\|_\phi )$. 
According to Lemma \ref{Le:HKE} and Markov property, for each $\mu\in S_{\phi}$ and $t\in (0,1]$, $\psi(\mu)_t$ admits a density $p^\mu_t(y)$ w.r.t.  Lebesgue measure and 
$$
p^\mu_t(y) = \int_{\R^d} p^\mu(0,x;t,y) \pi(\d x). 
$$
 So for each fixed $t\in (0,1]$,  the relatively compactness of $\{\psi(\mu)_t: \mu\in S_{\phi}\}$ in $(\cP_{\phi}(\R^d), \|\cdot\|_\phi )$ is equivalent to the relatively compactness of $\{\phi p^\mu_t: \mu\in S_{\phi}\}$  in $L^1(\R^d)$. By \eqref{Eq:TSE},\eqref{Eq:growth} and Lebesgue's  dominated convergence theorem, 
\begin{align*}
&\lim_{R\to\infty} \sup_{\mu\in S_\phi} \int_{|y|>R} \phi p^\mu_t (y) \d y\\
\overset{\eqref{Eq:TSE}}{\lesssim}&\lim_{R\to\infty} \int_{\R^d}\pi(\d x)\int_{\R^d}\varrho_\lambda(t,x-y)\1_{B_R^c}(y) \phi(y) \d y\\
=&\int_{\R^d}\pi(\d x) \lim_{R\to\infty}  \int_{|y|>R}\phi(x-y)\varrho_\lambda(t,y) \d y=0.
\end{align*}
On the other hand,  by  \eqref{Eq:Holder-y},  \eqref{Eq:growth} and Lebesgue's  dominated convergence theorem, we obtain that there is a constant $\gamma\in (0, \a\wedge\gamma_0)$ such that 
\begin{align*}
&\sup_{\mu\in S_\phi}\int_{\R^d} |\phi p^\mu_t(y+h)-\phi p^\mu_t(y)| \d y \\
\leq & \sup_{\mu\in S_\phi}\int_{\R^d} | p^\mu_t(y+h)- p^\mu_t(y)| \phi(y+h)\d y+ \sup_{\mu\in S_\phi}\int_{\R^d} p^\mu_t(y)| \phi(y+h)-\phi(y)| \d y \\
\overset{\eqref{Eq:Holder-y}}{\lesssim}& t^{-\frac{\gamma}{2}}|h|^\gamma \int_{\R^d} \pi(\d x)\int_{\R^d}\l[ \varrho_{\lambda}(t, x-y-h) +\varrho_{\lambda}(t, x-y)\r] \phi(y+h)\d y\\
&+  |h| \int_0^1\d \theta \int_{\R^d} \pi(\d x)\int_{\R^d} \varrho_{\lambda}(t,x-y) |\nabla \phi(y+\theta h)| \d y \\
\overset{\eqref{Eq:growth}}{\lesssim}& (t^{-\frac{\gamma}{2}}  |h|^{\gamma} +|h|)\<\phi, \pi\>  \lesssim t^{-\frac{\gamma}{2}}  |h|^{\gamma}. 
\end{align*}
Thanks to Fr\'echet-Kolmogorov's theorem, we get the desired result. 

\medskip 

(3): Assume $\mu^n$ is a sequence in $S_\phi$ and  $\mu^n\to \mu$ in $S_\phi$ as $n\to \infty$. By our assumption, we have 
$$
\eps_n:= \|a^{\mu^n}-a^\mu\|_{\cL^\infty_r}^{1-\eta}+ \|b^{\mu^n}-b^\mu\|_{\cL^p_q}\to 0\quad (n\to \infty). 
$$
Denote the heat kernel associated with $L^{\mu_n}=a^{\mu_n}_{ij}\p_{ij}+b^{\mu_n}_i\p_i$ and $L^\mu=a^{\mu}_{ij}\p_{ij}+b^{\mu}_i\p_i$ by $p^{\mu_n}(s, x;t,y)$ and $p^{\mu}(s,x;t,y)$, respectively. Then, for each $t\in (0,1]$
\begin{align*}
\|\psi(\mu^n)_t-\psi(\mu)_t\|_{\phi}=&\sup_{\|f\|_{L^\infty}\leq 1}\l|\bE f\phi(X_t^{\mu_n}) -\bE f\phi(X_{t}^{\mu})\r|\\
\leq& \int_{\R^d} \pi(\d x)\int_{\R^d}|p^{\mu_n}-p^{\mu}|(0,x; t, y)\phi(y) \d y\\
\overset{\eqref{Eq:stable}}{\lesssim}& \eps_n t^{-\frac{1}{r}}\int_{\R^d} \pi(\d x)\int_{\R^d}\phi(y) \varrho_\lambda(t,x-y) \d y\\
\overset{\eqref{Eq:growth}}{\lesssim}& \eps_n t^{-\frac{1}{r}} \<\phi, \pi\>  \to 0\quad (n\to\infty), 
\end{align*}
which implies $\psi$ is continuous from $S_{\phi}$ to $K$. 
\epf

\section{Uniqueness in law}\label{Sec:Uni}
In this section, we give our main result about the uniqueness of solutions to \eqref{Eq:MV}, which is addressed by comparing the densities of two solutions.  As we already mentioned in the introduction, the uniformly ellipticity of $a$ and the Lipschitz continuity of $\sigma$ in $m$ w.r.t. the total variation distance is not enough to ensure the uniqueness. Below, we illustrate this point with a simple example.
\bex[Nonuniqueness]\label{Ex-ex3}
Let $d=d_1=1$, $W_t$ be a standard Brownian motion on $\R$. Let 
$$
c_1=\frac{1}{\sqrt{2\pi}} \int_{-{2}}^{2} \e^{-\frac{|x|^2}{2}} \d x \approx 0.95, \ \ 
c_2=\frac{1}{\sqrt{8\pi}} \int_{-{2}}^{2} \e^{-\frac{|x|^2}{8}} \d x \approx 0.68
$$
and $(\lambda_1, \lambda_2)$ be the solution to the following linear system of  equations  
\begin{align*}
\l\{ \begin{aligned}
c_1\lambda_1+(1-c_1)\lambda_2=&1\\
c_2\lambda_1+(1-c_2)\lambda_2=&2
\end{aligned}
\r. 
\end{align*}
i.e. 
$$\lambda_2=\frac{2c_1-c_2}{c_1-c_2}>\lambda_1=\frac{2c_1-c_2-1}{c_1-c_2}>0. 
$$
Assume that 
$$
\Sigma(t, x)= \lambda_1 \1_{B_{2}}(x/\sqrt{t})+ \lambda_2 \1_{B_{2}^c}(x/\sqrt{t}), \quad \si(t, m):=\int_{\R} \Sigma(t,x) m(\d x) 
$$
and $b(t, x, m)\equiv 0$. According to our definition, $\sigma(t, m)$ satisfies $0 < \lambda_1 \leq \sigma(t, m) \leq \lambda_2 < \infty$, and the mapping from $m$ to $\sigma(t, m)$ exhibits uniform Lipschitz continuity concerning the total variation distance. Moreover, modulo an additive constant, $\Sigma(t, y) = \frac{\delta\sigma}{\delta m}(t, m)(y)$ remains uniformly bounded. However, through straightforward calculations, we find that
$$
\si(t, \mu_{W_t})=  \lambda_1 \mu_{W_1}(B_2)+ \lambda_2 \mu_{W_1}(B_2^c)=c_1\lambda_1+(1-c_1)\lambda_2=1
$$
and 
$$
\si(t, \mu_{2W_t})=  \lambda_1 \mu_{2W_1}(B_2)+ \lambda_2 \mu_{2W_1}(B_2^c)=c_2\lambda_1+(1-c_2)\lambda_2=2,  
$$
which imply that \eqref{Eq:MV} (with $\pi=\delta_0$) has at least two strong solutions: $W_t$ and $2W_t$. 
\eex

Inspired by the work of de Ranal-Frikha \cite{de2020strong, de2022well}, we establish our uniqueness result based on the assumption that $\frac{\delta a}{\delta m}(t, \cdot, m)(\cdot)$ exhibits uniform H\"older continuity for all $(t, m)$. The following theorem represents the main result of this section. 
\bt[Uniqueness]\label{Thm:Uni}
Let $(p,q)\in \sI_1$, $\alpha, \beta \in (0,1)$,  $\Lambda>1$ and $N_1, N_2>0$. Assume that 
\begin{enumerate}[(i)] 
\item  for all $m\in \cP_\phi(\R^d)$, 
\begin{align}\label{Aspt3}
\begin{aligned}
a(\cdot,\cdot,m)\in \mS(\Lambda, \a, N_1),  \quad \sup_{t\in [0,1]}\l\|\frac{\delta a}{\delta m} (t,\cdot, m)(\cdot)\r\|_{C^\beta(\R^{2d})}\leq N_1; 
\end{aligned}
\tag{\bf U$_\sigma$}
\end{align}
\item there is a nonnegative function $\ell \in L^	q([0,1];\R_+)$ such that for any $\mu \in S_\phi$ and $m, m'\in \cP_\phi(\R^d)$, 
\begin{align}\label{Aspt4}
\|b^\mu\|_{\cL^p_q}\leq N_2, \quad \|b(t,\cdot,m)-b(t,\cdot,m')\|_{\cL^p}\leq \ell(t)\|m-m'\|_{\phi} \tag{\bf U$_b$}. 
\end{align}
\end{enumerate}
Then for any $\pi\in\cP_\phi(\R^d)$, \eqref{Eq:MV} admits a unique weak solution. 
\et
\bpf
Note that \eqref{Aspt3} and \eqref{Aspt4} imply \eqref{Aspt1} and \eqref{Aspt2}, respectively, so by Lemma \ref{Le:Ex}, we only need to prove the uniqueness. 
Assume $X$ and $\widetilde X$ are two weak solutions to \eqref{Eq:MV}. Let $\mu_t=\mu_{X_t}$, $\widetilde \mu_t=\mu_{\widetilde X_t}$,  $L=a^\mu_{ij}\p_{ij}+b^\mu_{i}\p_i$, $\widetilde L=a^{\widetilde \mu}_{ij}\p_{ij}+b^{\widetilde \mu}_i\p_i$ and $p$, $\widetilde p$ be the heat kernels associated with $L$ and $\widetilde L$, respectively. As in the proof for Proposition \ref{Prop-Key}(i), we have 
\begin{align}\label{eq:rho-mu}
\begin{aligned}
&\sup_{t\in [0,1]} \<\phi, \mu_t\>= \sup_{t\in [0,1]}\int_{\R^{d}}\pi(\d x)  \int_{\R^d} p(0,x;t,y) \phi(y) \d y\\
\overset{\eqref{Eq:TSE}}{\lesssim }&  \sup_{t\in [0,1]}\int_{\R^d} \pi(\d x)\int_{\R^d}\varrho_\lambda(t,x-y) \phi(y) \d y 
\overset{\eqref{Eq:growth}}{\lesssim}   \<\phi, \pi\>. 
\end{aligned}
\end{align}

Assume $\mu_s=\widetilde \mu_s$ for some $s\in [0,1]$. Below we prove that there is a constant $T>0$ that only depends on $d, p ,q, \a, \beta, \Lambda, N_1, N_2, c$ and $\<\phi, \pi\>$ such that for any $t\in [s, 1\wedge (s+T)]$, $\mu_t=\widetilde \mu_t$. 

Define   
\begin{equation}\label{Eq:def-epsT}
q:=p-\widetilde p,\quad \eps(T):=\sup_{\substack{x,y\in \R^d;\\ t\in (s, 1\wedge(s+T)]}}\frac{|q|(s, x; t, y)}{\varrho_{\lambda}(t-s,x-y)}<\infty. 
\end{equation}
By the identity $p=p_0+p\otimes [(L-L_0)p_0]$ and the definition of $q$, one see that 
\begin{equation}\label{Eq:def-q}
\begin{aligned}
q=& (p_0-\widetilde p_0)+(p-\widetilde p) \otimes (L-L_0)p_0\\
&+ \widetilde p\otimes [(L- L_0)-(\widetilde L-\widetilde L_0)]p_0+ \widetilde  p\otimes (\widetilde L-\widetilde L_0)(p_0- \widetilde p_0)\\
=&: J_0+J_1+J_2+J_3. 
\end{aligned}
\end{equation}
The Markovian property, along with our assumption that $\mu_s = \widetilde \mu_s$ yields  
\begin{equation}\label{eq:mu-mu}
(\mu_t-\widetilde \mu_t)(\d y)=\int_{\R^d} q(s,x;t,y)\mu_{s}(\d x) \, \d y, \quad \forall t\in (s,1]. 
\end{equation}
To analyze $J_0$, we employ \eqref{Eq:def-epsT}-\eqref{eq:mu-mu}, and our assumptions outlined in \eqref{Aspt3}. This yields the following estimate for all $\tau$ within the interval $[s, 1\wedge(s+T)]$: 
\begin{equation*}
\begin{aligned}
&|a(\tau,y, \mu_\tau)-a(\tau,y,\widetilde \mu_\tau)|\\
\overset{\eqref{Eq:f1-f2}}{=}&\left| {\int_{0}^{1}\!\!\! \int_{\mathbb{R}^{d}}\frac{\delta a}{\delta m}\left(\tau, y, \lambda\mu_\tau+(1-\lambda)\widetilde \mu_\tau\right)\left(z\right)}\ (\mu_\tau-\widetilde \mu_\tau)(\d z)  \d \lambda\right|\\
=&\bigg| \int_{0}^{1}\!\!\! \int_{\mathbb{R}^{d}} \Big[\frac{\delta a}{\delta m}\left(\tau, y, \lambda\mu_\tau+(1-\lambda)\widetilde \mu_\tau\right)(z)-\frac{\delta a}{\delta m}\left(\tau, y, \lambda\mu_\tau+(1-\lambda)\widetilde \mu_\tau\right)(x)\Big] \\
& \qquad\qquad\qquad\qquad\qquad\qquad \qquad\qquad\qquad\qquad\qquad\qquad\qquad  \cdot(\mu_\tau-\widetilde \mu_\tau)(\d z)  \d \lambda\bigg|\\
{\leq} &{\int_{\R^d}\!\int_{0}^{1}\!\!\!\int_{\mathbb{R}^{d}}\left|\frac{\delta a}{\delta m}\left(\tau, y, \lambda\mu_\tau+(1-\lambda)\widetilde \mu_\tau\right)\left(z\right)-\frac{\delta a}{\delta m}\left(\tau, y, \lambda\mu_\tau+(1-\lambda)\widetilde \mu_\tau\right)\left(x\right)\right|}\\
&\qquad\qquad\qquad\qquad\qquad\qquad\qquad\qquad\qquad\qquad\qquad \cdot|q \left(s, x; \tau,  z\right)| \d z \  \d \lambda \ \mu_s(\d x)\\
{\lesssim}&  \eps(T)\int_{\R^d} \mu_s(\d x)\int_{\R^d}  |x-z|^\beta\varrho_\lambda(\tau-s, x-z)\d z\lesssim  \eps(T) (\tau-s)^{\frac{\beta}{2}}.  
\end{aligned}
\end{equation*}
Like the proof for \eqref{eq:p0-p0-1}, by the above estimate, we obtain that for any $ t \in[s,1\wedge(s+T)]$, 
\begin{equation}\label{eq:J0}
\begin{aligned}
&|J_0|(s,x;t,y)=|p_0-\widetilde p_0|(s,x;t,y)\\
\overset{\eqref{eq:p0-p0-1}}{\lesssim}&  \varrho_{2\lambda, -1}(t-s, x-y) \int_s^t |a(\tau,y, \mu_\tau)-a(\tau,y,\widetilde \mu_\tau)|\d \tau \\
\lesssim & \eps(T) \varrho_{2\lambda,-1}(t-s, x-y)  \int_s^t (\tau-s)^{\frac{\beta}{2}}\d\tau\\
\lesssim& \eps(T) T^{\frac{\beta}{2}}\varrho_{\lambda}(t-s, x-y) . 
\end{aligned}
\end{equation}
Similarly, for any $k=0,1,2$, we have 
\begin{equation}\label{eq:p0k-p0k}
|\nabla^k(p_0-\widetilde p_0) |(s,x;t,y)
\lesssim \eps(T) T^{\frac{\beta}{2}}\varrho_{\lambda, -k}(t-s, x-y) . 
\end{equation}
For $J_1$, suppose that $t\in [s,1\wedge(s+T)]$, then 
\begin{equation}\label{eq:J1}
\begin{split}
|J_1|(s,x;t,y) \overset{\eqref{Eq:kp0}}{\lesssim} & \int_s^t\!\!\! \int_{\R^d} |q (s,x;\tau, z)|\varrho_{2\lambda, \a-2}(t-\tau, z-y) \d z\d \tau \\
&+\int_s^t\!\!\! \int_{\R^d} |q(s,x;\tau, z)|\, |b(\tau,z)|\varrho_{2\lambda, -1}(t-\tau, z-y) \d z\d \tau\\
\overset{\eqref{Eq:def-epsT}}{\lesssim}& \eps(T) \int_s^t (t-\tau)^{\frac{\a}{2}-1} \d \tau \int_{\R^d} \varrho_\lambda(\tau-s, x-z) \varrho_{\lambda}(t-\tau, z-y)\d z\\
&+\eps(T) \int_s^t\!\!\!\int_{\R^d} \varrho_\lambda(\tau-s,x-z) |b(\tau,z)| \varrho_{2\lambda,-1}(t-\tau, z-y) \d z \d \tau \\
\overset{\eqref{Eq:kato}}{\lesssim} & \eps(T) \l(T^{\frac{\a}{2}}+K_{|b|}^1(T)\r)\varrho_{\lambda}(t-s; x-y). 
\end{split}
\end{equation}
For $J_2$. Noting that for any $\tau\in [s, 1]$, 
\begin{equation*}
\begin{aligned}
 &\l| [(L- L_0)-(\widetilde L-\widetilde L_0)] p_0\r|(\tau,z;t, y)\\
\leq & \l|\l[a_{ij}(\tau, z,\mu_\tau)-a_{ij}(\tau, y, \mu_\tau)]-[a_{ij}(\tau, z, \widetilde\mu_\tau)-a_{ij}(\tau, y, \widetilde\mu_\tau)\r]\r| |\p_{z_iz_j}p_0(\tau,z;t, y)|\\
&+ \l| b_i^{\mu}(\tau, z)- b_i^{\widetilde \mu}(\tau, z)\r| | \p_{z_i}p_0(\tau, z; t, y)|\\
\leq& \left| {\int_{0}^{1}\!\!\! \int_{\mathbb{R}^{d}}\l[\frac{\delta a}{\delta m}\left(\tau, z, \lambda\mu_\tau+(1-\lambda)\widetilde \mu_\tau\right)-\frac{\delta a}{\delta m}\left(\tau, y, \lambda\mu_\tau+(1-\lambda)\widetilde \mu_\tau\right) \r]\left(z'\right)}\r. \\
&\l. \cdot (\mu_\tau-\widetilde \mu_\tau)(\d z')  \, \d \lambda\right|\cdot\varrho_{2\lambda, -2}(t-\tau,z-y)\\
&+\l| b^{\mu}(\tau, z)- b^{\widetilde \mu}(\tau, z)\r| \varrho_{2\lambda, -1}(t-\tau,z-y) \\
\overset{\eqref{Aspt3},\eqref{eq:mu-mu}}{\lesssim} &\varrho_{2\lambda, -2}(t-\tau,z-y)  \int_{\R^{2d}} |y-z|^\beta |q|(s,x';\tau,z') \d z' \mu_s(\d x')\\
&+ \varrho_{2\lambda, -1}(t-\tau,z-y)\l| b^{\mu}(\tau, z)- b^{\widetilde \mu}(\tau, z)\r|, 
\end{aligned}
\end{equation*}
we have 
\begin{equation*}
\begin{aligned}
&|J_2|(s,x;t,y)\\
\overset{\eqref{Eq:def-epsT}}{\lesssim}& \eps(T)\int_s^t\!\!\!\int_{\R^{3d}}  \varrho_{\lambda}(\tau-s,x-z)  \cdot \varrho_{2\lambda, -2}(t-\tau,z-y)\\
 & \quad \quad \quad \quad \quad \quad \quad \quad |y-z|^\beta\varrho_{\lambda}(\tau-s, x'-z') \d z ' \mu_s(\d x') \d z \d\tau \\
&+ \int_s^t\!\!\!\int_{\R^d} \varrho_{\lambda}(\tau-s,x-z) \l| [b^{\mu}(\tau, z)- b^{\widetilde \mu}(\tau, z)]\1_{[s,t]}(\tau)\r| \varrho_{2\lambda, -1}(t-\tau,z-y) \d z\d\tau\\
\overset{\eqref{Eq:kato}}{\lesssim} & \eps(T) \int_s^t\!\!\!\int_{\R^d} \varrho_{\lambda}(\tau-s,x-z)\varrho_{\lambda, \beta-2}(t-\tau,z-y)  \d z\d\tau+K^1_{|b^{\mu}-b^{\widetilde \mu}|\1_{[s,t]}}(T)\varrho_\lambda(t-s, x-y) , \\
\lesssim& \eps(T) T^{\frac{\beta}{2}}  \varrho_{\lambda}(t-s, x-y)+ K^1_{|b^{\mu}-b^{\widetilde \mu|}\1_{[s,t]}}(T) \varrho_\lambda(t-s, x-y). 
\end{aligned}
\end{equation*} 
By \eqref{Eq:KvsL}, the last term of above inequalities can be controlled by 
\begin{align*}
&\om(T)\l\| (b^\mu-b^{\widetilde\mu})\1_{[s,t]}\r\|_{\cL^p_q}  \varrho_{\lambda}  \\
\overset{\eqref{Aspt4}}{\lesssim} & \om(T) \sup_{\tau\in [s,t]} \|\mu_\tau-\widetilde \mu_\tau\|_{\phi} \varrho_{\lambda}\\
\lesssim& \om(T) \l(\sup_{\tau\in[s,1\wedge(s+T)]}\int_{\R^d}  \mu_s(\d x) \int_{\R^{d}}|q(s,x;\tau,z)| \phi(z)\d z\r)\varrho_{\lambda}  \\
\overset{\eqref{Eq:growth}}{\lesssim} & \om(T)  \eps(T) \<\phi, \mu_s\> \varrho_{\lambda} \overset{\eqref{eq:rho-mu}}{\lesssim}\om(T) \eps(T) \<\phi, \pi\>\varrho_{\lambda} ,  
\end{align*}
where $\om(T)\to 0$ as $T\to 0$. Thus, for each  $t\in [s,1\wedge (s+T)]$
\begin{equation}\label{eq:J2}
|J_2|(s,x;t,y)\lesssim \l(T^{\frac{\beta}{2}}   + \om(T) \<\phi, \pi\> \r)\eps(T) \varrho_{\lambda}(t-s,x-y). 
\end{equation}
For $J_3$, by \eqref{eq:p0k-p0k} and \eqref{Eq:def-epsT}, we have 
\begin{align*}
&(\widetilde L-\widetilde L_0)(p_0-\widetilde p_0) (\tau,z;t,y)\\
\lesssim &\eps(T) T^{\frac{\beta}{2}} \l[ \varrho_{2\lambda, \a-2}(t-\tau, z-y)+  |b^{\widetilde \mu}(\tau,z)|\varrho_{2\lambda, -1}(t-\tau, z-y) \r]. 
\end{align*}
Thus, 
\begin{equation}\label{eq:J3}
|J_3|(s,x;t,y)\lesssim \eps(T) T^{\frac{\beta}{2}} \l(T^{\frac{\a}{2}}+K_{|b^{\widetilde \mu}|}^1(T)\r) \varrho_{\lambda}(t-s; x-y). 
\end{equation}
Combining \eqref{Eq:def-q}-\eqref{eq:J3}, we obtain 
\begin{align*}
\eps(T) =&\sup_{\substack{y\in \R^d;\\t\in [s, 1\wedge(s+T)] }}\frac{|q|(s, x; t, y)}{\varrho_{\lambda}(t-s,x-y)}\leq \sup_{\substack{y\in \R^d;\\t\in [s, 1\wedge s+T] }} \frac{\sum_{i=0}^{3}J_i(s,x;t,y)}{\varrho_{\lambda}(t-s,x-y)}\\
\leq& C\l(1+\<\phi, \pi\>\r)\om(T)\eps(T), 
\end{align*}
where $C$ is a constant that only depends on $d, p ,q, \a, \beta, \Lambda, N_1, N_2, c$. This  implies that there exists $T>0$, which is independent of $s$ and such that $\eps(T)\equiv 0$. Thus,  
$$
(\mu_t-\widetilde \mu_t)(\d y)= \int_{\R^d} q(s,x;t,y) \mu_s(\d x)\equiv 0, \quad \forall t\in [s,1\wedge (s+T)]. 
$$
Since $\mu_0=\widetilde \mu_0=\pi$, we obtain $\mu_t=\widetilde \mu_t$ for all $t\in [0,1]$. So we complete our proof. 
\epf

Our Proposition \ref{Prop:spf} presented in the introduction is a corollary of Proposition \ref{Thm:Uni}. 
\bpf[Proof of Theorem \ref{Prop:spf}]
Take $\phi=1$. Suppose $\si, b$ are given by \eqref{Eq:A-B}. By the assumptions on $\Sigma$, we can observe that
$$
\sup_{(t,m)\in [0,1]\times \cP(\R^d)} \l(\|a(t,\cdot, m)\|_{C^\a(\R^d)}+\l\|\frac{\delta a}{\delta m} (t,\cdot, m)(\cdot)\r\|_{C^\alpha(\R^{2d})} \r) \leq C \sup_{t\in [0,1]}\|\Sigma(t, \cdot)\|_{C^\a}<\infty, 
$$
which implies $a$ satisfies \eqref{Aspt3}. On the other hand, by the Minkowski inequality, one sees that 
\begin{align*}
\|b^\mu\|_{\cL^p_q}= \l\|\int_{\R^d}B(t, x-y) \mu_t(\d y)\r\|_{L^q_tL^p_x}\leq \|B\|_{L^q_tL^p_x}  
\end{align*}
and 
\begin{align*}
\|b(t,\cdot,m)-b(t,\cdot,m')\|_{\cL^p}\leq \|B(t,\cdot)\|_{L^p_x}\|m-m'\|_{{\rm TV}}.
\end{align*}
Thus, $b$ satisfies \eqref{Aspt4} with $\ell(t)=\|B(t,\cdot)\|_{L^p_x}\in L^q([0,1])$. So we obtain the desired result due to Theorem \ref{Thm:Uni}. 
\epf

Next we briefly comment on the strong well-posedness of \eqref{Eq:MV}. 

 Under the same conditions of Lemma \ref{Le:Ex}, by \cite[Theorem 1.1]{xia2020lqlp}, if in addition 
\begin{align}\label{Eq:Dsig}
\|\nabla_x\si^{\mu}\|_{\cL^p_q}<\infty, \quad \forall\mu\in S_\phi, 
\end{align}
then \eqref{Eq:MV} has at least one strong solution, provided that $\pi\in \cP_{\phi}(\R^d)$. Furthermore, under the conditions specified in Theorem \ref{Thm:Uni}, along with the constraints outlined in \eqref{Eq:Dsig}, \eqref{Eq:MV} admits a unique strong solution, provided that $\pi\in \cP_{\phi}(\R^d)$.  

Lastly, let us also mention that \eqref{Eq:MV} is closely related to the following nonlinear Fokker-Planck equation (NFPE): 
\begin{equation}\label{Eq:NFPE}
\p_t\mu_t-\p_{ij}(a_{ij}(t,x, \mu_t)\mu_t) \mu_t +\p_i (b_i(t,x, \mu_t)\mu_t) =0, \  \mu_0=\pi\in \cP_\phi(\R^d). 
\end{equation}
 Recall that a curve $\mu: [0,1]\to \cP(\R^d)$ is narrowly continuous, if for every $f\in C_b(\R^d)$, $t\mapsto \<f, \mu_t\>$ is continuous. The following conclusion regarding NFPE \eqref{Eq:NFPE} is an consequence drawn from \cite[Theorem 1.1]{xia2020lqlp}, \cite[Theorem 5.1]{rockner2021DDSDE} and our primary results .
 
\bc\label{Cor-NFPE}
\begin{enumerate}
\item Assume $a, b$ satisfy \eqref{Aspt1} and \eqref{Aspt2} (or \eqref{Eq:Esig} and \eqref{Eq:Eb}), then for any $\pi\in \cP_{\phi}(\R^d)$, \eqref{Eq:NFPE} has  at least one weak solution $\mu:[0,1]\to \cP(\R^d)$ such that $\mu$ is narrow continuous and satisfies the following Krylov's type estimate, 
\begin{equation}\label{Eq:Krylov}
\int_0^1\!\!\!\int_{\R^d} |f|(t, x) \mu_{t} (\d x)\d t\leq C \|f\|_{\cL^p_q}. 
\end{equation}
\item Assume $a, b$ satisfy \eqref{Aspt3} and \eqref{Aspt4}, then for each $\pi\in\cP_\phi(\R^d)$,  \eqref{Eq:NFPE} admits a unique weak solution in $C([0,1];\cP(\mR^d))$. 
\end{enumerate}
\ec

\appendix

\section{Proof for Lemma \ref{Le:LSDE} and Lemma \ref{Le:Kato}}
  \setcounter{equation}{0}
  \renewcommand\theequation{A.\arabic{equation}}

\begin{proof}[Proof of Lemma \ref{Le:LSDE}]
	The proof is essentially contained in \cite{xia2020lqlp}. 
	We provide an outline of the proofs below for the reader's convenience. To simplify, we can assume without loss of generality that $s=0$. Leveraging the given assumptions and Theorem 3.2 of \cite{xia2020lqlp}, for any  $c,f\in C_c^\infty(\mR^{d+1})$, there exists a unique solution $u\in \cH^{2,p}_q$ with $\partial_t u\in \cL^p_q$, which satisfies the following parabolic equation:
	$$
	\p_t v+\frac{1}{2}\si_{ik}\si_{jk}\p_{ij} v+ b_i\p_i v+c v+f, \quad v(1)=0. 
	$$
	By the generalized It\^o's formula (see \cite[Lemma 4.1]{xia2020lqlp}), we have 
	\begin{align*}
		&\d [\e^{\int_0^t c(s, X_s) \d s}v(t, X_t)]\\
		=& - \e^{\int_0^t c(s, X_s)\d s} f(t, X_t)+ \e^{\int_0^t c(s, X_s)\d s} \nabla_x v(t, X_t)\, \si(t, X_t)\, \d W_t. 
	\end{align*}
	This implies 
	$$
	\bE \int_{0}^{1} \e^{\int_0^t c(s, X_s) \d s} f\left(t, X_{t}\right) \mathrm{d} t= \bE v(0, X_0). 
	$$
	From this, we derive all the weak solutions have the same finite dimensional distributions. So the proof is complete.
\end{proof}

\medskip

\begin{proof}[Proof of Lemma \ref{Le:Kato}]
	(i) Noting 
	$$
	|x|^\gamma \varrho_{\lambda, 0}(t, x) =t^{(-d+\gamma)/2}   (|x|^2/t)^{\gamma/2} \e^{-\lambda |x|^2/t}, 
	$$
	by the elementary inequality $a \e^{-\lambda a} \geq C(\lambda ,\kappa) \e^{-\kappa\lambda a} \ (\forall a\geq 0)$, we get the desired estimate. 
	
	(ii) If $|x|<\sqrt{t}$, then 
	$$
	\varrho_{\lambda, -\beta}(t,x)\leq t^{-(d+\beta)/2}\lesssim \eta_\beta(t,x). 
	$$ 
	On the other hand, if $|x|>\sqrt{t}$, then $$
	\varrho_{\lambda,-\beta}(t,x)= |x|^{-d-\beta} \l[ \l(|x|/\sqrt{t}\r)^{d+\beta}\e^{-\lambda |x|^2/t} \r]\lesssim  |x|^{-d-\beta}\lesssim \eta_\beta(t,x).
	$$
	
	(iii) For any $k\geq 1$, there is an integer $N_k\asymp  k^{d-1}$ and a sequence of unit balls $\{B_1(x_{k, i})\}_{i=1}^{N_k}$ such that $B_{k+1}\backslash B_k\subseteq \bigcup_{i=1}^{N_k} B_1(x_{k, i})$. 
	For any $T\in (0,1)$ and $f\in \cL^p_q$ with $(p,q)\in \sI_\beta$,  by H\"older's inequality and the fact $(\beta+d/p)q/(2q-2)<1$,  we have 
	\begin{align*}
		&\int_0^T\!\!\!\int_{\R^d} |f(s, y)| \eta_\beta(s, y) \d y \d s \lesssim \int_0^Ts^{-(d+\beta)/2}\d s\int_{B_{\sqrt{s}}} |f(s, y)| \d y\\
		&\quad \quad \quad +  \int_0^T\!\!\!\int_{B_1\backslash B_{\sqrt s}} \frac{|f(s, y)|}{|y|^{d+\beta}} \d y \d s + \sum_{k=1}^\infty\int_0^T\!\!\!\int_{B_{k+1}\backslash B_k}\frac{|f(s, y)|}{|y|^{d+\beta}}  \d y\d s\\
		\lesssim & \int_0^T  \|f(s, \cdot )\1_{B_1}\|_{L^p}s^{-(\beta+d/p)/2} \d s+ \int_0^T  \|f(s, \cdot)\1_{B_1}\|_{L^p}\l(\int_{\sqrt s}^1 r^{\frac{d+\beta}{p-1}-\beta-1} \d r\r)^{\frac{p-1}{p}} \d s \\
		&+ \sum_{k=1}^\infty k^{-d-\beta} \sum_{i=1}^{N_k} \int_0^T\!\!\!\int_{\R^d} |f(s,y)| \1_{B_1(x_{k,i})}(y) \ \d y \d s \\
		\lesssim & \|f\|_{\cL^p_q(T)}\sqrt{T}^{2-\beta-\frac{d}{p}-\frac{2}{q}}+ \|f\|_{\cL^p_q(T)} T^{1-\frac{1}{q}} \sum_{k=1}^\infty k^{-1-\beta}\lesssim T^{\frac{1}{2}(2-\beta-\frac{d}{p}-\frac{2}{q})} \|f\|_{\cL^p_q(T)}. 
	\end{align*}
	
	(iv) 
	The essence of the proof for \eqref{Eq:kato} can be found in \cite[Lemma 3.1]{zhang1997gaussian}. We include the proof here for the reader's convenience, assuming, for the sake of simplicity, that $s=0$. Let 
	\begin{equation*}
		I=\int_{0}^{t}\!\!\!\int_{\R^d} \frac{\varrho_{\lambda}(\tau, x-z)}{\tau^{\beta'/2}}|b(\tau,z)| \frac{\varrho_{2\lambda}(t-\tau, z-y)}{(t-\tau)^{\beta/2}} \d z \d \tau. 
	\end{equation*}
	Then 
	\begin{align*}
		I=I_1+I_2:=& \int_{0}^{\vartheta t}\!\!\!\int_{\R^d}  \frac{\varrho_{\lambda}(\tau, x-z)}{\tau^{\beta'/2}}|b(\tau,z)| \frac{\varrho_{2\lambda}(t-\tau, z-y)}{(t-\tau)^{\beta/2}} \d z \d \tau\\
		&+\int_{\vartheta t}^t \int_{\R^d} \frac{\varrho_{\lambda}(\tau, x-z)}{\tau^{\beta'/2}}|b(\tau,z)| \frac{\varrho_{2\lambda}(t-\tau, z-y)}{(t-\tau)^{\beta/2}} \d z \d \tau, 
	\end{align*}
	where $\vartheta= \frac{3}{4}-\frac{1}{\sqrt{2}}$.  For $I_1$, we write 
	$$
	I_1= I_{11}+I_{12}:= \int_0^{\vartheta t}\!\!\!\int_{|z-y|>|x-y|/\sqrt{2}}\cdots+ \int_0^{\vartheta t}\!\!\!\int_{|z-y|\leq |x-y|/\sqrt{2}}\cdots. 
	$$
	When $|z-y|> |x-y|/\sqrt{2}$ and $\tau\in [0,\vartheta t]$, we have 
	$$
	{\varrho_{2\lambda}(t-\tau, z-y)} \leq \varrho_{\lambda}(t-\tau, x-y)\lesssim \varrho_\lambda(t,x-y). 
	$$ 
	This yields  
	\begin{align*}
		I_{11}\lesssim& \varrho_{\lambda, -\beta'}(t,x-y)\int_{0}^{\vartheta t}\!\!\!\int_{\R^d}  \frac{\varrho_{\lambda}(\tau, x-z)}{\tau^{\frac{\beta'}{2}}(t-\tau)^{\frac{\beta-\beta'}{2}}}|b(\tau,z)|  \d z\d \tau\\
		\lesssim &\varrho_{\lambda, -\beta'}(t,x-y)\int_{0}^{t}\!\!\!\int _{\R^d} \frac{\varrho_{\lambda}(\tau, x-z)}{\tau^{\beta/2}}|b(\tau, z)|  \d z\d \tau\\
		\overset{\eqref{Eq:rhoVeta}}{\lesssim}  &  K_{|b|}^\beta(t)\varrho_{\lambda, -\beta'}(t,x-y).  
	\end{align*}
	When $|z-y|\leq |x-y| /\sqrt{2}$ and $\tau\in[0,\vartheta t]$, we have 
	$$
	|x-z|^2/2\tau\geq (1-1/\sqrt{2})^2\l(\tfrac{3}{2}-\sqrt{2}\r)^{-1} |x-y|^2/t=|x-y|^2/t, 
	$$ 
	which implies 
	$$
	\varrho_{\lambda}(\tau, x-z)\leq \varrho_{\lambda/2} (\tau, x-z) \, \e^{-|x-y|^2/t}. 
	$$
	Thus, 
	\begin{align*}
		I_{12}\lesssim&  \e^{-|x-y|^2/t} \int_{0}^{\vartheta t}\!\!\!\int_{\R^d}  \frac{\varrho_{\lambda/2}(\tau, x-z)}{\tau^{\beta'/2}} |b(\tau,z)| \frac{1}{(t-\tau)^{(d+\beta)/2}} \d z \d \tau \\
		\lesssim& \varrho_{\lambda, -\beta'} (t, x-y) \int_{0}^{t}\!\!\!\int _{\R^d} |b(\tau,z)| \frac{\varrho_{\lambda/2}(\tau, x-z)}{\tau^{\beta/2}} \d z \d \tau \\
		\overset{\eqref{Eq:rhoVeta}}{\lesssim}  & K_{|b|}^\beta(t)\varrho_{\lambda, -\beta'} (t, x-y). 
	\end{align*}
	For $I_2$, by the elementary inequality,
	$$
	\frac{a^2}{\a}+ \frac{b^2}{1-\a} \geq (a+b)^2, \quad \forall a, b\geq 0, \a\in (0,1)
	$$
	we obtain 
	$$
	\frac{|x-z|^{2}}{\tau}+\frac{|z-y|^{2}}{t-\tau} \geq \frac{|x-y|^{2}}{t}, \quad 0<\tau<t. 
	$$
	Combining the above inequality and the fact that $\tau\geq \vartheta t$, one sees that 
	$$
	\frac{\varrho_{\lambda}(\tau, x-z)}{\tau^{\beta'/2}}\cdot\frac{\varrho_{2\lambda}(t-\tau, z-y)}{(t-\tau)^{\beta/2}}\lesssim \frac{\varrho_{\lambda}(t, x-y)}{t^{\beta'/2}} \frac{\varrho_\lambda(t-\tau, z-y)}{(t-\tau)^{\beta/2}}. 
	$$
	Hence, 
	\begin{align*}
		I_2 \lesssim& \varrho_{\lambda, -\beta'}(t,x-y)\int_{0}^t \!\!\!\int_{\R^d} |b(t-\tau,z)| \frac{\varrho_{\lambda}(\tau, y-z) }{\tau^{\beta/2}} \d z \d \tau\\
		\overset{\eqref{Eq:rhoVeta}}{\lesssim} & K^{\beta}_{|b|}(t)\varrho_{\lambda, -\beta'}(t,x-y)
	\end{align*}
\end{proof}

\section*{Acknowledgements}
The author would like to thank Professor Michael R\"ockner , Professor Xicheng Zhang and Xianliang Zhao for many helpful discussions and  comments.

The author finished the main part of the article while working at Bielefeld University and received funding from the German Research Foundation (DFG) through the Collaborative Research Centre(CRC) 1283 ``Taming uncertainty and profiting from randomness and low regularity in analysis, stochastics and their applications" during that time.

\bibliographystyle{alpha}
\bibliography{myref}

\begin{thebibliography}{CHXZ17}

\bibitem[BP18]{bauer2018strong}
Thilo Meyer-Brandis Bauer, Martin and Frank Proske.
\newblock Strong solutions of mean-field stochastic differential equations with
  irregular drift.
\newblock {\em Electronic Communications in Probability}, 23(132):1--35, 2018.

\bibitem[BR18]{barbu2018probabilistic}
Viorel Barbu and Michael R{\"o}ckner.
\newblock Probabilistic representation for solutions to nonlinear
  fokker--planck equations.
\newblock {\em SIAM Journal on Mathematical Analysis}, 50(4):4246--4260, 2018.

\bibitem[BR20]{barbu2020nonlinear}
Viorel Barbu and Michael R\"{o}ckner.
\newblock From nonlinear {F}okker-{P}lanck equations to solutions of
  distribution dependent {SDE}.
\newblock {\em Ann. Probab.}, 48(4):1902--1920, 2020.

\bibitem[BR21]{barbu2021uniqueness}
Viorel Barbu and Michael R\"{o}ckner.
\newblock Uniqueness for nonlinear {F}okker-{P}lanck equations and weak
  uniqueness for {M}c{K}ean-{V}lasov {SDE}s.
\newblock {\em Stoch. Partial Differ. Equ. Anal. Comput.}, 9(3):702--713, 2021.

\bibitem[CDLL19]{cardaliaguet2015master}
Pierre Cardaliaguet, Fran\c{c}ois Delarue, Jean-Michel Lasry, and Pierre-Louis
  Lions.
\newblock {\em The master equation and the convergence problem in mean field
  games}, volume 201 of {\em Annals of Mathematics Studies}.
\newblock Princeton University Press, Princeton, NJ, 2019.

\bibitem[CdRF22]{de2022well}
Paul-Eric Chaudru~de Raynal and Noufel Frikha.
\newblock Well-posedness for some non-linear {SDE}s and related {PDE} on the
  {W}asserstein space.
\newblock {\em J. Math. Pures Appl. (9)}, 159:1--167, 2022.

\bibitem[CHXZ17]{chen2017heat}
Zhen-Qing Chen, Eryan Hu, Longjie Xie, and Xicheng Zhang.
\newblock Heat kernels for non-symmetric diffusion operators with jumps.
\newblock {\em Journal of Differential Equations}, 263(10):6576--6634, 2017.

\bibitem[dR20]{de2020strong}
Paul-Eric~Chaudru de~Raynal.
\newblock Strong well posedness of {M}c{K}ean-{V}lasov stochastic differential
  equations with {H}{\"o}lder drift.
\newblock {\em Stochastic Processes and their Applications}, 130(1):79--107,
  2020.

\bibitem[FG17]{guilin2017kac}
Nicolas Fournier and Arnaud Guillin.
\newblock From a {K}ac-like particle system to the {L}andau equation for hard
  potentials and {M}axwell molecules.
\newblock {\em Ann. Sci. \'{E}c. Norm. Sup\'{e}r. (4)}, 50(1):157--199, 2017.

\bibitem[{F}ig08]{figalli2008existence}
{A}lessio {F}igalli.
\newblock {E}xistence and uniqueness of martingale solutions for {S}{D}{E}s
  with rough or degenerate coefficients.
\newblock {\em {J}ournal of {F}unctional {A}nalysis}, 254(1):109--153, 2008.

\bibitem[Fri08]{friedman2008partial}
Avner Friedman.
\newblock {\em Partial differential equations of parabolic type}.
\newblock Courier Dover Publications, 2008.

\bibitem[Fun84]{funaki1984certain}
Tadahisa Funaki.
\newblock A certain class of diffusion processes associated with nonlinear
  parabolic equations.
\newblock {\em Zeitschrift f{\"u}r Wahrscheinlichkeitstheorie und Verwandte
  Gebiete}, 67(3):331--348, 1984.

\bibitem[HW19]{huang2019distribution}
Xing Huang and Feng-Yu Wang.
\newblock Distribution dependent {{SDE}}s with singular coefficients.
\newblock {\em Stochastic Processes and their Applications},
  129(11):4747--4770, 2019.

\bibitem[Jou97]{jourdain1997diffusions}
Benjamin Jourdain.
\newblock Diffusions with a nonlinear irregular drift coefficient and
  probabilistic interpretation of generalized {B}urgers' equations.
\newblock {\em Probability and Statistics}, 1(1):339--355, 1997.

\bibitem[Lac18]{lacker2018strong}
Daniel Lacker.
\newblock On a strong form of propagation of chaos for {M}ckean-{V}lasov
  equations.
\newblock {\em Electronic Communications in Probability}, 23(45):1--11, 2018.

\bibitem[LM16]{li2016weak}
Juan Li and Hui Min.
\newblock Weak solutions of mean-field stochastic differential equations and
  application to zero-sum stochastic differential games.
\newblock {\em SIAM Journal on Control and Optimization}, 54(3):1826--1858,
  2016.

\bibitem[MRS15]{manita2015uniqueness}
Oxana~A Manita, Maxim~S Romanov, and Stanislav~V Shaposhnikov.
\newblock On uniqueness of solutions to nonlinear fokker--planck--kolmogorov
  equations.
\newblock {\em Nonlinear Analysis}, 128:199--226, 2015.

\bibitem[MV20]{mishura2020existence}
Yuliya Mishura and Alexander Veretennikov.
\newblock Existence and uniqueness theorems for solutions of
  {M}c{K}ean-{V}lasov stochastic equations.
\newblock {\em Theory Probab. Math. Statist.}, (103):59--101, 2020.

\bibitem[Pat19]{pata2019fixed}
Vittorino Pata.
\newblock {\em Fixed point theorems and applications}, volume 116.
\newblock Springer, 2019.

\bibitem[Rud73]{rudin1973functional}
Walter Rudin.
\newblock {\em Functional analysis}.
\newblock McGraw-hill, 1973.

\bibitem[RZ21]{rockner2021DDSDE}
Michael R\"{o}ckner and Xicheng Zhang.
\newblock Well-posedness of distribution dependent {SDE}s with singular drifts.
\newblock {\em Bernoulli}, 27(2):1131--1158, 2021.

\bibitem[ST85]{shiga1985central}
Tokuzo Shiga and Hiroshi Tanaka.
\newblock Central limit theorem for a system of {M}arkovian particles with mean
  field interactions.
\newblock {\em Zeitschrift f{\"u}r Wahrscheinlichkeitstheorie und Verwandte
  Gebiete}, 69(3):439--459, 1985.

\bibitem[Szn91]{sznitman1991topics}
Alain-Sol Sznitman.
\newblock Topics in propagation of chaos.
\newblock In {\em Ecole d'{\'e}t{\'e} de probabilit{\'e}s de Saint-Flour
  XIX—1989}, pages 165--251. Springer, 1991.

\bibitem[{T}re12]{trevisan2012well}
{D}ario {T}revisan.
\newblock {W}ell-posedness of multidimensional diffusion processes with weakly
  differentiable coefficients.
\newblock {\em {E}lectronic {J}ournal of {P}robability}, 21, 2012.

\bibitem[Vil08]{villani2008optimal}
C{\'e}dric Villani.
\newblock {\em Optimal transport: old and new}, volume 338.
\newblock Springer Science \& Business Media, 2008.

\bibitem[Wan18]{wang2018distribution}
Feng-Yu Wang.
\newblock Distribution dependent {{SDE}}s for {L}andau type equations.
\newblock {\em Stochastic Processes and their Applications}, 128(2):595--621,
  2018.

\bibitem[XXZZ20]{xia2020lqlp}
Pengcheng Xia, Longjie Xie, Xicheng Zhang, and Guohuan Zhao.
\newblock ${L}^q({L}^p)$-theory of stochastic differential equations.
\newblock {\em Stochastic Processes and their Applications}, 130(8):5188--5211,
  2020.

\bibitem[Zha97]{zhang1997gaussian}
Qi~S Zhang.
\newblock Gaussian bounds for the fundamental solutions of $\nabla ({A}\nabla
  u)+ {B}\nabla u- u_t= 0$.
\newblock {\em manuscripta mathematica}, 93(1):381--390, 1997.

\bibitem[ZZ18]{zhang2018singular}
Xicheng Zhang and Guohuan Zhao.
\newblock {S}ingular {B}rownian diffusion processes.
\newblock {\em {C}ommunications in {M}athematics and {S}tatistics},
  6(4):533--581, 2018.

\bibitem[ZZ21]{zhang2020stochastic}
Xicheng Zhang and Guohuan Zhao.
\newblock Stochastic {L}agrangian path for {L}eray’s solutions of 3{D}
  {N}avier--{S}tokes equations.
\newblock {\em Communications in Mathematical Physics}, 381(2):491--525, 2021.

\end{thebibliography}

\end{document}